\documentclass[11pt,a4paper,fullpage]{article}
\usepackage[english]{babel}
\usepackage[utf8]{inputenc}
\usepackage{fancyhdr}
\usepackage[numbers]{natbib}

\usepackage{graphicx}
\usepackage{amsmath}
\usepackage{extsizes}
\usepackage{relsize}
\usepackage{multicol}
\usepackage{ragged2e}
\usepackage{fontenc}
\usepackage{mathcomp}
\usepackage{latexsym}
\usepackage{amssymb}
\usepackage{times}
\newtheorem{Teorema}{Theorem}[section]

\newtheorem{Definicija}[Teorema]{Definition}
\newtheorem{Posledica}[Teorema]{Corollary}

\newtheorem{Lema}[Teorema]{Lemma}
\newtheorem{Primedba}[Teorema]{Remark}
\usepackage{dsfont}
\newcommand\scalemath[2]{\scalebox{#1}{\mbox{\ensuremath{\displaystyle #2}}}} 
\numberwithin{equation}{section}
\hyphenchar\font=-1
\begin{document}
	\begin{center}
	{\LARGE Intersections of various spectra of operator matrices\par}	
\vspace{1cm}
	{\large Nikola Sarajlija\footnote{e-mail address: nikola.sarajlija@dmi.uns.ac.rs\\University of Novi Sad, Faculty of Sciences, Novi Sad, 21000, Serbia}}


\vspace{5mm}
	{\large \today\par}
\end{center}
\vspace{1cm}
\begin{abstract}
This paper is concerned with general $n\times n$ upper triangular operator matrices with given diagonal entries. We characterize perturbations of the left (right) essential spectrum, the essential spectrum, as well as the left (right) the Weyl spectrum of such operators, thus generalizing and improving some results of \cite{CAO}, \cite{CAO2}, \cite{SUNTHEORY}, \cite{SUN}, \cite{LIandDU}, \cite{ZHANG} from $n=2$ to case of $n\geq2$. We show that related results hold in arbitrary Banach spaces without assuming separability, thus extending and improving recent results from \cite{WU}, \cite{WU2}. For the lack of adjointness and orthogonality in Banach spaces, we use an adequate alternative: the concept of inner generalized inverse. 
\end{abstract}
\textbf{Keywords: }essential spectrum, Weyl spectrum, perturbation theory, $n\times n$, regular operators
\section{Introduction}

In the last few decades considerable attention has been devoted to the study of spectral properties of operator matrices, having in mind their governing importance in various areas of mathematics (see for example \cite[Chapter VIII]{DUNFORD}). One has soon realized that one way for successful work on problems arising in spectral theory, is to see operator matrices as entries of simpler blocks. Block operator matrices, and especially upper triangular operator matrices, have been extensively studied by numerous authors (see \cite{CAO}, \cite{CAO2}, \cite{SUNTHEORY}, \cite{DU}, \cite{HAN}, \cite{LI}, \cite{KINEZI}, \cite{WU}, \cite{WU2} and many others...). The reason for this lies in the fact that if an operator $T$ is acting on a direct sum of Banach spaces, it takes the upper triangular form under condition that certain number of those spaces is invariant for $T$.

Development of this topic began in the last century, and is of great importance ever since. In the beginning, authors have only considered the case of $2\times2$ operator matrices. Pioneering work in that direction was the article of Du and Jin from 1994 (\cite{DU}) treating the usual spectrum. Han et al. have generalized their result to Banach spaces (\cite{HAN}), and Lee has proved some facts concerning the Weyl spectrum (\cite{LI}). Afterwards, Djordjevi\'{c} in 2002 gave some characterizations for $2\times2$ upper triangular operators to be Fredholm, Weyl, and Browder (\cite{SUNTHEORY}). After that, many authors have explored various properties of $2\times2$ block operators in a connection with intersection of spectra, Weyl and Browder type theorems, etc. (see for example \cite{KINEZI}, \cite{CAO}).

Investigation of spectral properties of general $n\times n$ operators began no sooner than 2013, when Zguitti published the article dealing with some properties of the Drazin spectrum (\cite{ZGUITTI}). Huang et al. continued his work in 2016 by investigating properties of the point, residual, and continuous spectrum of $n\times n$ matrix operators (\cite{HUANG}).  Fredholm and the Weyl spectrum of such operators have been studied by Wu and Huang in \cite{WU}, \cite{WU2} only a few years ago, and this article is concerned with generalizing their results from separable Hilbert to arbitrary Banach spaces.

Partial operator matrix is an operator matrix with some entries specified and the others unknown. Most commonly, those entries are linear and bounded operators with appropriate domains. We mention that there has been some interest in block operators with unbounded entries lately, see \cite{NEOGRANICENI3}, \cite{NEOGRANICENI2}, \cite{NEOGRANICENI1}, but we shall not pursue this point any further. Usually, it is required that unknown entries satisfy some special properties, and finding those entries is the task known as the completion of the matrix. It is sometimes a rather difficult task, but there are some ways to make this task easier. Perturbation of an upper triangular operator matrix is an intersection of its certain spectrum, some of the entries being fixed, for example the diagonal ones, and the others chosen arbitrarly with appropriate domains. Finding perturbation is a very useful pre-step to the completion of a matrix. This article is ought to be just a such pre-step.

Article is organized as follows. In Section 2 we give some basic concepts and terminology that we use in the rest of this work. We introduce notion of embedded spaces that serves us as a surrogate to a lack of separability, and we present some of its variants. Perturbation results related to Fredholmness and Weylness of upper triangular block operators belong to Section 3 and Section 4, respectively. We supply each of those sections with an additional subsection that contains an answer to the Question 3 raised in Section 2.

\section{Notation and preliminaries}

Let $X, Y, X_1,...,X_n$ be arbitrary Banach spaces. We use the notation $\mathcal{B}(X_i,X_j)$ for the collection of all linear and bounded operators from $X_i$ to $X_j$. Particularly, $\mathcal{B}(X_i)=\mathcal{B}(X_i,X_i)$, and $X'=\mathcal{B}(X,\mathds{C})$ is the dual space of $X$. If $T\in\mathcal{B}(X_i,X_j)$, then by $\mathcal{N}(T)$ and $\mathcal{R}(T)$ we denote the kernel and range of $T$. Those sets are subspaces of $X_i$ and $X_j$, respectively, and $\mathcal{N}(T)$ is closed.

If $T\in\mathcal{B}(X,Y)$, then its dual operator is $T'\in\mathcal{B}(Y',X')$ defined by $T'f(x):=f(Tx),\ f\in Y'$. We shall use some properties of dual operators on several occasions. For example, it is known that $\Vert T\Vert=\Vert T'\Vert$ and $T\mapsto T'$ is an isometric isomorphism of $\mathcal{B}(X,Y)$ into $\mathcal{B}(X',Y')$. For other basic features of dual operators we recommend \cite{DUNFORD}.

For $U\subseteq X$ we define set $U^\circ\subseteq X'$, and for $V\subseteq X'$ we define set $^\circ V\subseteq X$ as
$$U^\circ:=\lbrace f\in X':\ f\upharpoonleft_U=0\rbrace.$$
$$^\circ V:=\lbrace x\in X:\ V(x)=0\rbrace.$$
$U^\circ$ and $^\circ V$ are called the left and right annihilator of $U$ and $V$, respectively. Above all interesting features that hold for annihilators, we point out only one needed in the proof of Lemma \ref{VELIKALEMA}, see \cite{SETER}:
\begin{equation}\label{ANULATORI}
\overline{\mathcal{R}(A)}=^\circ[\mathcal{N}(A')],
\end{equation}
where $A\in\mathcal{B}(X,Y)$.

Let $D_1\in\mathcal{B}(X_1),\ D_2\in\mathcal{B}(X_2),...,D_n\in\mathcal{B}(X_n)$ be given. We denote by $T_n^d(A)$ an $n\times n$ partial upper triangular operator matrix of the form
\begin{equation}\label{OSNOVNA}
T_n^d(A)=
\begin{bmatrix} 
    D_1 & A_{12} & A_{13} & ... & A_{1,n-1} & A_{1n}\\
    0 & D_2 & A_{23} & ... & A_{2,n-1} & A_{2n}\\
    0 &  0 & D_3 & ... & A_{3,n-1} & A_{3n}\\
    \vdots & \vdots & \vdots & \ddots & \vdots & \vdots\\
    0 & 0 & 0 & ... & D_{n-1} & A_{n-1,n}\\
    0 & 0 & 0 & ... & 0 & D_n      
\end{bmatrix}\in\mathcal{B}(X_1\oplus X_2\oplus\cdots\oplus X_n),
\end{equation}
where $A:=(A_{12},\ A_{13},...,\ A_{ij},...,\ A_{n-1,n})$ is an operator tuple consisting of unknown variables $A_{ij}\in\mathcal{B}(X_j,X_i)$, $1\leq i<j\leq n,\ n\geq2$. For convenience, we denote by $\mathcal{B}_n$ the collection of all such tuples $A$. One easily verifies that for $T_n^d(A)$ of the form (\ref{OSNOVNA}), its dual operator matrix $T_n^d(A)'$ is given by
\begin{equation}\label{DUALNA}
T_n^d(A)'=
\scalemath{0.9}{
\begin{bmatrix} 
    D_1' & 0 & 0 & ... & 0 & 0\\
    A_{12}' & D_2' & 0 & ... & 0 & 0\\
    A_{13}' &  A_{23}' & D_3' & ... & 0 & 0\\
    \vdots & \vdots & \vdots & \ddots & \vdots & \vdots\\
    A_{1,n-1}' & A_{2,n-1}' & A_{3,n-1}' & ... & D_{n-1}' & 0\\
    A_{1n}' & A_{2n}' & A_{3n}' & ... & A_{n-1,n}' & D_n'      
\end{bmatrix}}\in\mathcal{B}(X_1'\oplus X_2'\oplus\cdots\oplus X_n').
\end{equation}

There are several questions that we are interested in:

\noindent\textbf{Question 1.} Can we find an appropriate characterization for Fredholmness, Weylness, etc. for $T_n^d(A)$, in terms of Fredholmness, Weylness, etc. of its diagonal entries $D_i$?

\noindent\textbf{Question 2.} What can we say about $\bigcap\limits_{A\in\mathcal{B}_n}\sigma_*(T_n^d(A))$, when $\sigma_*$ runs through set $\lbrace\sigma_{le},\sigma_{re},\sigma_e,\sigma_{lw},\sigma_{rw}\rbrace$?

\noindent\textbf{Question 3.} Under what conditions the equality $\sigma_*(T_n^d(A))=\bigcup\limits_{i=1}^n\sigma_*(D_i)$ holds, $\sigma_*\in\lbrace\sigma_{le},\sigma_{re},\sigma_e, \sigma_{lw}, \sigma_{rw} \rbrace$?

Here we give complete answers to Questions 2 i 3, while Question 1 remains opened.  

We use $\mathcal{G}_l(X)$ and $\mathcal{G}_r(X)$ to denote the sets of all left and right invertible operators on $X$, respectively. It is well known that if $T\in\mathcal{B}(X)$, than $T\in\mathcal{G}_l(X)$ if and only if $\mathcal{N}(T)=\lbrace0\rbrace$ and $\mathcal{R}(T)$ is closed and complemented. Similarly, $T\in\mathcal{G}_r(X)$ if and only if $\mathcal{R}(T)=X$ and $\mathcal{N}(T)$ is complemented. The set of all invertible operators on $X$ is denoted by $\mathcal{G}(X)=\mathcal{G}_l(X)\cap\mathcal{G}_r(X)$.

Now we list some elementary notions from Fredholm theory (see \cite{ZANA}). Let $T\in\mathcal{B}(X)$, $\alpha(T)=\dim\mathcal{N}(T)$ and $\beta(T)=\dim X/\mathcal{R}(T)$. Quantities $\alpha$ and $\beta$ are called the nullity and deficiency of $T$, respectively, and in the case where at least one of them is finite we define $\mathrm{ind}(T)=\alpha(T)-\beta(T)$ to be the index of $T$. Notice that $\mathrm{ind}(T)$ may be $\pm\infty$ or integer. Families of left and right Fredholm operators, respectively, are defined as
$$
\begin{aligned}
\Phi_l(X)=\lbrace T\in\mathcal{B}(X):\ \alpha(T)<\infty\ and\ \mathcal{R}(T)\ is \ closed\rbrace
\end{aligned}
$$
and
$$
\begin{aligned}
\Phi_r(X)=\lbrace T\in\mathcal{B}(X): \beta(T)<\infty\rbrace.
\end{aligned}
$$
The set of Fredholm operators is
$$\Phi(X)=\Phi_l(X)\cap\Phi_r(X)=\lbrace T\in\mathcal{B}(X): \alpha(T)<\infty\ and\ \beta(T)<\infty\rbrace.$$
Families of left and right Weyl operators, respectively, are defined as
$$\Phi_l^-(X)=\lbrace T\in\Phi_l(X):\ \mathrm{ind}(T)\leq0\rbrace$$
and
$$\Phi_r^+(X)=\lbrace T\in\Phi_r(X):\ \mathrm{ind}(T)\geq0\rbrace.$$
The set of Weyl operators is 
$$\Phi_0(X)=\Phi_l^-(X)\cap\Phi_r^+(X)=\lbrace T\in\Phi(X):\ \mathrm{ind}(T)=0\rbrace.$$

Corresponding spectra of an operator $T\in\mathcal{B}(X)$ are defined as follows:\\
-the left spectrum: $\sigma_{l}(T)=\lbrace\lambda\in\mathds{C}: \lambda-T\not\in\mathcal{G}_{l}(X)\rbrace$;\\
-the right spectrum: $\sigma_{r}(T)=\lbrace\lambda\in\mathds{C}: \lambda-T\not\in\mathcal{G}_{r}(X)\rbrace$;\\
-the spectrum: $\sigma(T)=\lbrace\lambda\in\mathds{C}: \lambda-T\not\in\mathcal{G}(X)\rbrace$;\\
-the left essential spectrum: $\sigma_{le}(T)=\lbrace\lambda\in\mathds{C}: \lambda-T\not\in\Phi_{l}(X)\rbrace$;\\
-the right essential spectrum: $\sigma_{re}(T)=\lbrace\lambda\in\mathds{C}: \lambda-T\not\in\Phi_{r}(X)\rbrace$;\\
-the essential spectrum: $\sigma_{e}(T)=\lbrace\lambda\in\mathds{C}: \lambda-T\not\in\Phi(X)\rbrace$;\\
-the left the Weyl spectrum: $\sigma_{lw}(T)=\lbrace\lambda\in\mathds{C}: \lambda-T\not\in\Phi_l^-(X)\rbrace$;\\
-the right the Weyl spectrum: $\sigma_{rw}(T)=\lbrace\lambda\in\mathds{C}: \lambda-T\not\in\Phi_r^+(X)\rbrace$;\\
-the Weyl spectrum: $\sigma_{w}(T)=\lbrace\lambda\in\mathds{C}: \lambda-T\not\in\Phi_0(X)\rbrace$.\\
All of these spectra are compact nonempty subsets of the complex plane.\\
We write $\rho_l(T),\rho_r(T), \rho(T), \rho_{le}(T), \rho_{re}(T), \rho_e(T), \rho_{lw}(T), \rho_{rw}(T), \rho_{w}(T)$ for the corresponding complements of the sets above, respectively.

The following lemma imposes a connection between $T$ and its dual operator $T'$ in terms of nullity and deficiency of $T$. This claim will be crucial at some points.

\begin{Lema}(\cite[p. 7-8]{CARADUS})\label{VEZA}
For $T\in\mathcal{B}(X)$ with closed range the following holds:\\
$(a)$ $\alpha(T)=\beta(T'),\beta(T)=\alpha(T')$;\\
$(b)$ $T\in\Phi_l(X)$ if and only if $T'\in\Phi_r(X')$;\\
$(c)$ $T\in\Phi_r(X)$ if and only if $T'\in\Phi_l(X')$;\\
$(d)$ $\mathrm{ind}(T')=-\mathrm{ind}(T).$
\end{Lema}

We emphasize the fact that our results are ought to hold in arbitrary Banach spaces, not just the separable ones. In order to prove main theorems which concern perturbation of various spectra of $T_n^d(A)$, we introduce a concept that will compensate loss of separability: the notion of embedded spaces. To our knowledge, this condition was first used in this context by D. S. Djordjevi\'{c} in 2002 (see \cite[Theorem 4.4]{SUNTHEORY}).

\begin{Definicija}\label{UTAPANJE}
Let $X$ and $Y$ be Banach spaces. We say that $X$ can be embedded in $Y$ and write $X\preceq Y$ if and only if there exists a left invertible operator $J:X\rightarrow Y.$
\end{Definicija}
\begin{Primedba}
Obviously, $X\preceq Y$ if and only if there exists a right invertible operator $J_1:Y\rightarrow X$.\\
If $H$ and $K$ are Hilbert spaces, than $H\preceq K$ if and only if $\dim H\leq \dim K$. Here, $\dim H$ stands for the orthogonal dimension of $H$.
\end{Primedba}

For the left (right) the Weyl spectrum we shall need another variant of the definition above (see \cite[Definition 4.2]{SUNTHEORY}).
\begin{Definicija}\label{ESENCIJALNOUTAPANJE}
Let $X$ and $Y$ be Banach spaces. We say that $X$ can be essentially embedded in $Y$ and write $X\prec Y$ if and only if:\\
(a) $X\preceq Y$;\\
(b) for every $T\in\mathcal{B}(X,Y)$, $Y/\mathcal{R}(T)$ is infinite dimensional.
\end{Definicija}
\begin{Primedba}
If $H$ and $K$ are Hilbert spaces, than $H\prec K$ if and only if $\dim H<\dim K$ and $K$ is infinite dimensional, where $\dim H$ is the orthogonal dimension of $H$.
\end{Primedba}

For the essential spectrum we need another form of Definition \ref{UTAPANJE}, that we introduce right now.
\begin{Definicija}\label{JAKOUTAPANJE}
Let $X$ and $Y$ be Banach spaces. We say that $X$ can be strongly embedded in $Y$ and write $X\preceq_s Y$ if and only if:\\
(a) $X\preceq Y$;\\
(b) there exists $J\in\mathcal{B}(X,Y)$ such that $Y/\mathcal{R}(J)$ is of finite dimension.
\end{Definicija}
\begin{Primedba}\label{SLICNISU}
Notice the similarity with the notion of isomorphism up to a finite dimensional subspace introduced in \cite{SUNTHEORY}. Correct link between these two is the following: if Banach spaces $X$ and $Y$ are isomorphic up to a finite dimensional subspace, then at least one of them can be strongly embedded into another.
\end{Primedba}
If $X$ and $Y$ are isomorphic up to a finite dimensional subspace, we write $$X\cong Y\ (u.f.d.s)$$ for brevity.

One important difference between Hilbert and Banach spaces is that closed subspace of a Hilbert space is always complemented $(\mathcal{H}=M\oplus M^\bot)$. This is not true for the case of Banach spaces. Since we would like to prove our results by decomposing Banach spaces in an appropriate way, we introduce the following condition.
\begin{Definicija}\label{KOMPLEMENTI}
We say that $D\in\mathcal{B}(X)$ satisfy the complements condition if and only if:\\
$(a)$ $\mathcal{N}(D)$ is complemented;\\
$(b)$ $\mathcal{R}(D)$ is closed and complemented.
\end{Definicija}

We immediately give an example of one important class of operators which satisfies the complements condition.

We say that an operator $T\in\mathcal{B}(X_i,X_j)$ is (relatively) regular if and only if there is $T'\in\mathcal{B}(X_j,X_i)$ such that $TT'T=T$ holds. In that case we say $T'$ is inner generalized inverse of $T$. Notice that existence of such $T'$ does not imply its uniqueness. One can prove the following characterization (see \cite[Corollary 1.1.5]{DJR}):
\begin{Teorema}
$T\in\mathcal{B}(X_i,X_j)$ is regular if and only if $\mathcal{N}(T)$ and $\mathcal{R}(T)$ are closed and complemented subspaces of $X_i$ and $X_j$, respectively.
\end{Teorema}
Thus, it follows that regular operators satisfy the complements condition introduced in Definition \ref{KOMPLEMENTI}.\\
It is important to highlight that operators in sets $\mathcal{G}_l(X),\ \mathcal{G}_r(X),\ \Phi_l(X),\ \Phi_r(X)$ are regular operators.\\
Moreover, it is easily proved that, following upper terminology, $TT'$ and $T'T$ are both projections, and so we have decompositions (see \cite[Theorem 1.1.3]{DJR})
\begin{equation}\label{REGULARNISU}X_i=\mathcal{N}(T)\oplus\mathcal{R}(T'T),\quad X_j=\mathcal{N}(TT')\oplus\mathcal{R}(T).\end{equation}

We continue with some introductory lemmas.
\begin{Lema}\label{VELIKALEMA}
Let $D_1\in\mathcal{B}(X_1),\ D_2\in\mathcal{B}(X_2),...,D_n\in\mathcal{B}(X_n)$. Then
$$\overline{\mathcal{R}(T_n^d(A)-\lambda)}=X_1\oplus X_2\oplus\cdots\oplus X_n\quad for\  every\  A\in\mathcal{B}_n$$
if and only if
$$\overline{\mathcal{R}(D_1-\lambda)}=X_1,\ \overline{\mathcal{R}(D_2-\lambda)}=X_2,...,\overline{\mathcal{R}(D_n-\lambda)}=X_n.$$
\end{Lema}
\textbf{Proof. }Since necessity holds obviously, we prove only sufficiency. Hence, assume that $\overline{\mathcal{R}(D_1-\lambda)}=X_1,\ \overline{\mathcal{R}(D_2-\lambda)}=X_2,...,\overline{\mathcal{R}(D_n-\lambda)}=X_n.$ Regarding the relation (\ref{ANULATORI}) we get
$\mathcal{N}(D_1'-\lambda)=\mathbf{0},\  \mathcal{N}(D_2'-\lambda)=\mathbf{0},...,\ \mathcal{N}(D_n'-\lambda)=\mathbf{0}.$ But then $\mathcal{N}(T_n^d(A)'-\lambda)\subseteq\bigcup\limits_{i=1}^{n}\mathcal{N}(D_i'-\lambda)$ implying $\mathcal{N}(T_n^d(A)'-\lambda)=\mathbf{0}$, and again with regards to (\ref{ANULATORI}) we have $\overline{\mathcal{R}(T_n^d(A)-\lambda)}=X_1\oplus X_2\oplus\cdots\oplus X_n$. $\square$
\begin{Lema}\label{POMOCNALEMA}
Let $T_n^d(A)\in\mathcal{B}(X_1\oplus\cdots\oplus X_n).$ Then:
\begin{itemize}
\item[(i)] $\sigma_{le}(D_1)\subseteq\sigma_{le}(T_n^d(A))\subseteq\bigcup\limits_{k=1}^n\sigma_{le}(D_k)$;
\item[(ii)] $\sigma_{re}(D_n)\subseteq\sigma_{re}(T_n^d(A))\subseteq\bigcup\limits_{k=1}^n\sigma_{re}(D_k);$
\item[(iii)]$\sigma_{lw}(T_n^d(A))\subseteq\bigcup\limits_{k=1}^n\sigma_{lw}(D_k)$;
\item[(iv)]$\sigma_{rw}(T_n^d(A))\subseteq\bigcup\limits_{k=1}^n\sigma_{rw}(D_k)$.
\end{itemize}
\end{Lema}

We continue with two more lemmas that will be essential in proving theorems to come.
\begin{Lema}(\cite[section 16, Lemma 2]{MULLER})\label{PRVALEMA} Let $T,F\in\mathcal{B}(X_1,X_2)$. If $F$ is a finite rank operator, then $\mathcal{R}(T+F)$ is closed if and only if $\mathcal{R}(T)$ is closed.
\end{Lema}
\begin{Lema}(\cite[Lemma 1.3]{APOSTOL} and \cite[Theorem 5.2]{KATO})\label{DRUGALEMA}
Let $D_i\in\mathcal{B}(X_i)$, $i=1,2$ and $A_{12}\in\mathcal{B}(X_2,X_1)$ with $\overline{\mathcal{R}(D_1)}=X_1$. If $\mathcal{R}(T_2^d(A))$ is closed, then $\mathcal{R}(D_2)$ is closed.
\end{Lema}

\section{Perturbations of left (right) essential spectrum and essential spectrum of $T_n^d(A)$}
Hereafter, if $D\in\mathcal{B}(X_i,X_j)$ satisfies the complements condition, complements of $\mathcal{N}(T)$ and $\mathcal{R}(T)$ will be denoted by $'$ sign, that is $X_i=\mathcal{N}(T)\oplus\mathcal{N}(T)'$ and $X_j=\mathcal{R}(T)\oplus\mathcal{R}(T)'$. As already announced, our intention is to exploit the property (\ref{REGULARNISU}) of regular operators as a substitution for orthogonality in Hilbert spaces. However, we will formulate our main theorems in a wider context of Definition \ref{KOMPLEMENTI}, and then obtain corresponding results expressed in the language of (\ref{REGULARNISU}) as consequences.

We start with a result which characterizes perturbations of left essential spectrum of $T_n^d(A)$. 
\begin{Teorema}\label{LEVIFREDHOLM}
Let $D_1\in\mathcal{B}(X_1),\ D_2\in\mathcal{B}(X_2),...,D_n\in\mathcal{B}(X_n)$ with all $D_i,\ i=1,...,n $ satisfying the complements condition. Additionally, suppose that we have 
\begin{equation}\label{PRVA}
\mathcal{N}(D_k)\preceq X_j/\mathcal{R}(D_j)\quad for\  all\  n\geq k>j\geq 1.
\end{equation}
Then
$$
\bigcap\limits_{A\in\mathcal{B}_n}\sigma_{le}(T_n^d(A))=\sigma_{le}(D_1)\cup\Big(\bigcup\limits_{k=2}^n\Delta_k\Big),
$$
where
$$
\Delta_k:=\Big\lbrace\lambda\in\sigma_{le}(D_k):\  \sum\limits_{s=1}^{k-1}\beta(D_s-\lambda)<\infty\Big\rbrace,\ 2\leq k\leq n.
$$
\end{Teorema}
\textbf{Proof. }Inclusion $(\supseteq)$ is already proved in \cite{WU}. For convenience of a reader we provide the full proof, with a slight difference in that now we have different space decompositions.

Write $\Delta_1:=\sigma_{le}(D_1)$,
$
\Delta_{k1}=\lbrace\lambda\in\mathds{C}: \mathcal{R}(D_k-\lambda)\  is\  not\  closed\ and\\ \ \sum\limits_{s=1}^{k-1}\beta(D_s-\lambda)<\infty\rbrace
$
and
$
\Delta_{k2}=\lbrace\lambda\in\mathds{C}: \mathcal{R}(D_k-\lambda)\  is\  closed\  and\\  \sum\limits_{s=1}^{k-1}\beta(D_s-\lambda)<\infty,\ \alpha(D_k-\lambda)=\infty\rbrace$
for $2\leq k\leq n$. \\Lemma \ref{POMOCNALEMA} implies $\Delta_1\subseteq\bigcap\limits_{A\in\mathcal{B}_n}\sigma_{le}(T_n^d(A))$ . Let $\lambda\in\bigcup\limits_{k=2}^n\Delta_k$.\\
Let $\lambda\in\Delta_{k1}$, that is $\mathcal{R}(D_k-\lambda)$  is not closed and $\sum\limits_{s=1}^{k-1}\beta(D_s-\lambda)<\infty$, $k\in\lbrace2,...,n\rbrace$. Then for each $A\in\mathcal{B}_k$, operator matrix $T_k^d(A)$ as an operator from $X_1\oplus\mathcal{N}(D_2)'\oplus\mathcal{N}(D_2)\oplus\mathcal{N}(D_3)'\oplus\mathcal{N}(D_3)\oplus\cdots\oplus\mathcal{N}(D_k)'\oplus\mathcal{N}(D_k)$ into $\mathcal{R}(D_1)\oplus\mathcal{R}(D_1)'\oplus\mathcal{R}(D_2)\oplus\mathcal{R}(D_2)'\oplus\cdots\oplus\mathcal{R}(D_{k-1})\oplus\mathcal{R}(D_{k-1}')\oplus\mathcal{R}(D_k)\oplus\mathcal{R}(D_k)'$ admits the following block representation
\begin{equation}\label{MATRICA}
T_k^d(A)=\scalemath{0.85}
{\begin{bmatrix} 
    D_1^{(1)} & A_{12}^{(1)} & A_{12}^{(2)} & A_{13}^{(1)} & A_{13}^{(2)} & ... & A_{1,k-1}^{(1)} & A_{1,k-1}^{(2)} & A_{1k}^{(1)} & A_{1k}^{(2)}\\
    0 & A_{12}^{(3)} & A_{12}^{(4)} & A_{13}^{(3)} & A_{13}^{(4)} & ... & A_{1,k-1}^{(3)} & A_{1,k-1}^{(4)} & A_{1k}^{(3)} & A_{1k}^{(4)}\\
    0 & D_2^{(1)} & 0 & A_{23}^{(1)} & A_{23}^{(2)} & ... & A_{2,k-1}^{(1)} & A_{2,k-1}^{(2)} & A_{2k}^{(1)} & A_{2k}^{(2)}\\
    0 & 0 & 0 & A_{23}^{(3)} & A_{23}^{(4)} & ... & A_{2,k-1}^{(3)} & A_{2,k-1}^{(4)} & A_{2k}^{(3)} & A_{2k}^{(4)}\\
    0 & 0 & 0 & D_{3}^{(1)} & 0 & ... & A_{3,k-1}^{(1)} & A_{3,k-1}^{(2)} & A_{3k}^{(1)} & A_{3k}^{(2)}\\
    0 & 0 & 0 & 0 & 0 & ... & A_{3,k-1}^{(3)} & A_{3,k-1}^{(4)} & A_{3k}^{(3)} & A_{3k}^{(4)}\\
    \vdots & \vdots & \vdots & \vdots & \vdots & \ddots & \vdots & \vdots & \vdots & \vdots\\
    0 & 0 & 0 & 0 & 0 & ... & D_{k-1}^{(1)} & 0 & A_{k-1,k}^{(1)} & A_{k-1,k}^{(2)}\\
    0 & 0 & 0 & 0 & 0 & ... & 0 & 0 & A_{k-1,k}^{(3)} & A_{k-1,k}^{(4)}\\
    0 & 0 & 0 & 0 & 0 & ... & 0 & 0 & D_k^{(1)} & 0\\
    0 & 0 & 0 & 0 & 0 & ... & 0 & 0 & 0 & 0\\
\end{bmatrix}}
\end{equation}
Assume that $\mathcal{R}(T_k^d(A)-\lambda)$ is closed. Then, according to Lemma \ref{PRVALEMA}, range of the following operator matrix
$$
\scalemath{0.85}{\begin{bmatrix}
    D_1^{(1)}-\lambda & A_{12}^{(1)} & A_{12}^{(2)} & A_{13}^{(1)} & A_{13}^{(2)} & ... & A_{1,k-1}^{(1)} & A_{1,k-1}^{(2)} & A_{1k}^{(1)} & A_{1k}^{(2)}\\
    0              & 0                   & 0                   & 0                   & 0                   & ... & 0                       & 0                       & 0                  & 0                   \\ 
    0 & D_2^{(1)}-\lambda & 0 & A_{23}^{(1)} & A_{23}^{(2)} & ... & A_{2,k-1}^{(1)} & A_{2,k-1}^{(2)} & A_{2k}^{(1)} & A_{2k}^{(2)}\\
    0 & 0 & 0 & 0 & 0 & ... & 0 & 0 & 0 & 0\\
    0 & 0 & 0 & D_{3}^{(1)}-\lambda & 0 & ... & A_{3,k-1}^{(1)} & A_{3,k-1}^{(2)} & A_{3k}^{(1)} & A_{3k}^{(2)}\\
    0 & 0 & 0 & 0 & 0 & ... & 0 & 0 & 0 & 0\\
    \vdots & \vdots & \vdots & \vdots & \vdots & \ddots & \vdots & \vdots & \vdots & \vdots\\
    0 & 0 & 0 & 0 & 0 & ... & D_{k-1}^{(1)}-\lambda & 0 & A_{k-1,k}^{(1)} & A_{k-1,k}^{(2)}\\
    0 & 0 & 0 & 0 & 0 & ... & 0 & 0 & 0 & 0\\
    0 & 0 & 0 & 0 & 0 & ... & 0 & 0 & D_k^{(1)}-\lambda & 0\\
    0 & 0 & 0 & 0 & 0 & ... & 0 & 0 & 0 & 0\\
\end{bmatrix}}
$$
is closed due to $\sum\limits_{s=1}^{k-1}\beta(D_s-\lambda)<\infty$, and hence range of the following operator matrix 
$$
\scalemath{0.85}{\begin{bmatrix}
D_1^{(1)}-\lambda & A_{12}^{(1)} & A_{12}^{(2)} & A_{13}^{(1)} & A_{13}^{(2)} & ... & A_{1,k-1}^{(1)} & A_{1,k-1}^{(2)} & A_{1k}^{(1)} & A_{1k}^{(2)}\\
    0 & D_2^{(1)}-\lambda & 0 & A_{23}^{(1)} & A_{23}^{(2)} & ... & A_{2,k-1}^{(1)} & A_{2,k-1}^{(2)} & A_{2k}^{(1)} & A_{2k}^{(2)}\\
    0 & 0 & 0 & D_{3}^{(1)}-\lambda & 0 & ... & A_{3,k-1}^{(1)} & A_{3,k-1}^{(2)} & A_{3k}^{(1)} & A_{3k}^{(2)}\\
    \vdots & \vdots & \vdots & \vdots & \vdots & \ddots & \vdots & \vdots & \vdots & \vdots\\
    0 & 0 & 0 & 0 & 0 & ... & D_{k-1}^{(1)}-\lambda & 0 & A_{k-1,k}^{(1)} & A_{k-1,k}^{(2)}\\
    0 & 0 & 0 & 0 & 0 & ... & 0 & 0 & D_k^{(1)}-\lambda & 0\\
\end{bmatrix}}
$$
is closed. Since $\overline{\mathcal{R}(D_s^{(1)}-\lambda)}=\overline{\mathcal{R}(D_s-\lambda)}$, $s=1,2,...,k-1$
, then
$$
\scalemath{0.85}{\mathcal{R}\left(\begin{bmatrix}D_1^{(1)}-\lambda & A_{12}^{(1)} & A_{12}^{(2)} & A_{13}^{(1)} & A_{13}^{(2)} & ... & A_{1,k-1}^{(1)} & A_{1,k-1}^{(2)}\\
    0 & D_2^{(1)}-\lambda & 0 & A_{23}^{(1)} & A_{23}^{(2)} & ... & A_{2,k-1}^{(1)} & A_{2,k-1}^{(2)}\\
    0 & 0 & 0 & D_{3}^{(1)}-\lambda & 0 & ... & A_{3,k-1}^{(1)} & A_{3,k-1}^{(2)}\\
    \vdots & \vdots & \vdots & \vdots & \vdots & \ddots & \vdots & \vdots\\
    0 & 0 & 0 & 0 & 0 & ... & D_{k-1}^{(1)}-\lambda & 0\\\end{bmatrix}\right)}
$$
is dense in $\overline{\mathcal{R}(D_1-\lambda)}\oplus\overline{\mathcal{R}(D_2-\lambda)}\oplus\cdots\oplus\overline{\mathcal{R}(D_{k-1}-\lambda)}$ from Lemma \ref{VELIKALEMA}. By Lemma \ref{DRUGALEMA}, it follows that $\mathcal{R}(D_k^{(1)}-\lambda)$ is closed, which is contradiction to the fact that $\mathcal{R}(D_k-\lambda)$ is not closed. This implies that $\mathcal{R}(T_k^d(A)-\lambda)$ is not closed for every $A\in\mathcal{B}_k$, and hence $\lambda\in\bigcap\limits_{A\in\mathcal{B}_k}\sigma_{le}(T_k^d(A))$. This proves $\bigcup\limits_{k=2}^n\Delta_{k1}\subseteq\bigcap\limits_{A\in\mathcal{B}_n}\sigma_{le}(T_n^d(A))$.

Let $\lambda\in\Delta_{k2}$, $k\in\lbrace 2,...,n\rbrace.$ Without loss of generality, put $\lambda=0$. We only consider $0\in\Delta_{k2}\setminus(\bigcup\limits_{s=2}^k\Delta_{s1}\bigcup\limits\Delta_1)$. Therefore, $0\in\rho_{le}(D_1)$, $\mathcal{R}(D_s)$ is closed for $s=2,...,k$, $\sum\limits_{s=1}^{k-1}\beta(D_s)<\infty$ and $\alpha(D_k)=\infty$. Then for each $A\in\mathcal{B}_k$, operator matrix $T_k^d(A)$ admits block representation as in (\ref{MATRICA}). Evidently, $D_1^{(1)},\ D_2^{(1)},...,D_k^{(1)}$ as in (\ref{MATRICA}) are all invertible. Thus, there exist invertible operator matrices $U$ and $V$ such that 
\begin{equation}\label{MATRICA2}
UT_k^d(A)V=\scalemath{0.85}{\begin{bmatrix}
D_1^{(1)} & 0 & 0 & 0 & 0 & ... & 0 & 0 & 0 & 0\\
0 & 0 & A_{12}^{(4)} & 0 & A_{13}^{(4)} & ... & 0 & A_{1,k-1}^{(4)} & 0 & A_{1k}^{(4)}\\
0 & D_2^{(1)} & 0 & 0 & 0 & ... & 0 & 0 & 0 & 0\\
0 & 0 & 0 & 0 & A_{23}^{(4)} & ... & 0 & A_{2,k-1}^{(4)} & 0 & A_{2k}^{(4)}\\
0 & 0 & 0 & D_3^{(1)} & 0 & ... & 0 & 0 & 0 & 0\\
0 & 0 & 0 & 0 & 0 & ... & 0 & A_{3,k-1}^{(4)} & 0 & A_{3k}^{(4)}\\
\vdots & \vdots & \vdots & \vdots & \vdots & \ddots & \vdots & \vdots & \vdots & \vdots\\
0 & 0 & 0 & 0 & 0 & ... & D_{k-1}^{(4)} & 0 & 0 & 0\\
0 & 0 & 0 & 0 & 0 & ... & 0 & 0 & 0 & A_{k-1,k}^{(4)}\\
0 & 0 & 0 & 0 & 0 & ... & 0 & 0 & D_k^{(1)} & 0\\
0 & 0 & 0 & 0 & 0 & ... & 0 & 0 & 0 & 0\\
\end{bmatrix}
}
\end{equation}
Note that $A_{ij}^{(4)}$ in (\ref{MATRICA2}) are not the original ones in (\ref{MATRICA}) in general, but we still use them for the convenience. Combine by (\ref{MATRICA2}), $T_k^d(A)$ is left Fredholm if and only if
\begin{equation}\label{MATRICA3}
\begin{bmatrix}
A_{12}^{(4)} & A_{13}^{(4)} & A_{14}^{(4)} & ... & A_{1k}^{(4)}\\
0          & A_{23}^{(4)} & A_{24}^{(4)} & ... & A_{2k}^{(4)}\\
0           & 0         & A_{34}^{(4)} & ... & A_{3k}^{(4)}\\
\vdots   &   \vdots & \vdots & \ddots & \vdots\\
0          &   0        &    0      & ... & A_{k-1,k}^{(4)}
\end{bmatrix}
:
\begin{bmatrix}
\mathcal{N}(D_2)\\
\mathcal{N}(D_3)\\
\mathcal{N}(D_4)\\
\vdots\\
\mathcal{N}(D_k)\\
\end{bmatrix}
\rightarrow
\begin{bmatrix}
\mathcal{R}(D_1)'\\
\mathcal{R}(D_2)'\\
\mathcal{R}(D_3)'\\
\vdots\\
\mathcal{R}(D_{k-1})'\\
\end{bmatrix}
\end{equation}
is left Fredholm. Since $\sum\limits_{s=1}^{k-1}\beta(D_s)<\infty$ and $\alpha(D_k)=\infty$, it follows that 
$$
\alpha\left(\begin{bmatrix}
A_{1k}^{(4)}\\
A_{2k}^{(4)}\\
A_{3k}^{(4)}\\
\vdots\\
A_{k-1,k}^{(4)}
\end{bmatrix}\right)=\infty,
$$
and hence operator defined in (\ref{MATRICA3}) is not left Fredholm for every $A\in\mathcal{B}_k$. This proves $\bigcup\limits_{k=2}^n\Delta_{k2}\subseteq\bigcap\limits_{A\in\mathcal{B}_n}\sigma_{le}(T_n^d(A))$.\\[3mm]

For the opposite inclusion, it suffices to prove that $0\not\in\bigcup\limits_{k=1}^n\Delta_k$ implies $0\not\in\bigcap\limits_{A\in\mathcal{B}_n}\sigma_{le}(T_n^d(A))$, i.e. there exists an operator tuple $A\in\mathcal{B}_n$ such that $T_n^d(A)$ is left Fredholm. Let $0\not\in\bigcup\limits_{k=1}^n\Delta_k$. Then, direct calculation show that $0\in\bigcup\limits_{k=1}^n\Omega_k$, where\\
$
\Omega_1=\lbrace\lambda\in\rho_{le}(D_1):\ \beta(D_1-\lambda)=\infty\rbrace\\
\Omega_k=\lbrace\lambda\in\bigcap\limits_{s=1}^k\rho_{le}(D_s):\ \beta(D_k-\lambda)=\infty\rbrace,\ 2\leq k\leq n-1,\\
\Omega_n=\bigcap\limits_{s=1}^n\rho_{le}(D_s).\\
$There are three cases to consider.\\[3mm]
\textbf{Case 1}: $0\in\Omega_1$\\

This means that $\alpha(D_1)<\infty$ ($D_1$ is left Fredholm), and $\beta(D_1)=\infty$. We find $A\in\mathcal{B}_n$ such that $\alpha(T_n^d(A))<\infty$ and $\mathcal{R}(T_n^d(A))$ is closed.  We choose $A=(A_{ij})_{1\leq i<j\leq n}$ so that $A_{ij}=0$ if $j-i\neq 1$, that is we place all nonzero operators of tuple $A$ on the superdiagonal. It remains to define $A_{ij}$ for $j-i=1$, $1\leq i<j\leq n$. First notice that $A_{i,i+1}:X_{i+1}\rightarrow X_i$. Since all of diagonal entries satisfy the complements condition, we know that $X_{i+1}=\mathcal{N}(D_{i+1})\oplus\mathcal{N}(D_{i+1})'$, $X_i=\mathcal{R}(D_i)'\oplus\mathcal{R}(D_i)$, and from assumption (\ref{PRVA}) we get $\mathcal{N}(D_{i+1})\preceq\mathcal{R}(D_i)'$. It follows that there is left invertible operator $J_{i}:\mathcal{N}(D_{i+1})\rightarrow\mathcal{R}(D_i)'$. We put $A_{i,i+1}=\begin{bmatrix}J_{i} & 0\\ 0 & 0\end{bmatrix}:\begin{bmatrix}\mathcal{N}(D_{i+1})\\ \mathcal{N}(D_{i+1})'\end{bmatrix}\rightarrow\begin{bmatrix}\mathcal{R}(D_i)'\\ \mathcal{R}(D_i)\end{bmatrix}$, and we implement this procedure for all $1\leq i\leq n-1$. Notice that $\mathcal{R}(D_i)$ is complemented to $\mathcal{R}(A_{i,i+1})$ for each $i=1,...,n-1.$

Now we have chosen our $A$, we show that $\mathcal{N}(T_n^d(A))\cong\mathcal{N}(D_1)$, implying $\alpha(T_n^d(A))=\alpha(D_1)<\infty$. Let us put $T_n^d(A)x=0$, where $x=x_1+\cdots+x_n\in X_1\oplus\cdots\oplus X_n$. The previous equality is then equivalent to the following system of equations
$$
\begin{bmatrix}D_1x_1+A_{12}x_2\\ D_2x_2+A_{23}x_3\\ \vdots\\ D_{n-1}x_{n-1}+A_{n-1,n}x_n\\ D_nx_n\end{bmatrix}=\begin{bmatrix}0\\ 0\\ \vdots\\ 0\\ 0\end{bmatrix}.
$$ 
The last equation gives $x_n\in\mathcal{N}(D_n)$. Since $\mathcal{R}(D_s)$ is complemented to $\mathcal{R}(A_{s,s+1})$ for all $1\leq s\leq n-1$, we have $D_sx_s=A_{s,s+1}x_{x+1}=0$ for all $1\leq s\leq n-1$. That is, $x_i\in\mathcal{N}(D_i)$ for every $i=1,...,n$, and $J_sx_{s+1}=0$ for every $1\leq s\leq n-1.$ Due to left invertibility of $J_s$ we get $x_s=0$ for $s=2,..., n$, which proves the claim. Therefore, $\alpha(T_n^d(A))=\alpha(D_1)<\infty$.

Secondly, we show that $\mathcal{R}(T_n^d(A))$ is closed. It is not hard to see that \\$\mathcal{R}(T_n^d(A))=\mathcal{R}(D_1)\oplus\mathcal{R}(J_{1})\oplus\mathcal{R}(D_2)\oplus\mathcal{R}(J_{2})\oplus\cdots\oplus\mathcal{R}(D_{n-1})\oplus\mathcal{R}(J_{n-1})\oplus\mathcal{R}(D_n)$. Furthermore, due to left invertibility of $J_i$'s, there exist closed subspaces $U_i$ of $\mathcal{R}(D_i)'$ such that $\mathcal{R}(D_i)'=\mathcal{R}(J_i)\oplus U_i$, $i=1,...,n -1$. It means that $X_1\oplus X_2\oplus\cdots\oplus X_n=\mathcal{R}(D_1)\oplus\mathcal{R}(J_{1})\oplus U_1\oplus\mathcal{R}(D_2)\oplus\mathcal{R}(J_{2})\oplus U_2\oplus\cdots\oplus\mathcal{R}(D_{n-1})\oplus\mathcal{R}(J_{n-1})\oplus U_{n-1}\oplus\mathcal{R}(D_n)$. Comparing these equalities, we conclude that $\mathcal{R}(T_n^d(A))$ is closed.\\

\noindent\textbf{Case 2}: $0\in\Omega_k$, $2\leq k\leq n-1$\\

This means that $\alpha(D_1),...,\alpha(D_k)<\infty$ ($D_1,..., D_k$ are left Fredholm), and $\beta(D_k)=\infty$. We find $A\in\mathcal{B}_n$ such that $\alpha(T_n^d(A))<\infty$ and $\mathcal{R}(T_n^d(A))$ is closed.  We choose $A=(A_{ij})_{1\leq i<j\leq n}$ analogously as in Case 1, with the only exception of choosing first $k-1$ operators on the superdiagonal to be zero. That is, $A_{i,i+1}=\begin{bmatrix}J_{i} & 0\\ 0 & 0\end{bmatrix}:\begin{bmatrix}\mathcal{N}(D_{i+1})\\ \mathcal{N}(D_{i+1})'\end{bmatrix}\rightarrow\begin{bmatrix}\mathcal{R}(D_i)'\\ \mathcal{R}(D_i)\end{bmatrix}$, $k\leq i\leq n$, and $A_{ij}=0$ otherwise, $1\leq i<j\leq n$, following notation of Case 1.

Now we have chosen our $A$, we show that $\mathcal{N}(T_n^d(A))\cong\mathcal{N}(D_1)\oplus\cdots\oplus\mathcal{N}(D_k)$. Let us put $T_n^d(A)x=0$, where $x=x_1+\cdots+x_n\in X_1\oplus\cdots\oplus X_n$. The previous equality is then equivalent to 
$$
\begin{bmatrix}D_1x_1\\ \vdots\\ D_{k-1}x_{k-1}\\ D_{k}x_{k}+A_{k,k+1}x_{k+1}\\ \vdots\\ D_{n-1}x_{n-1}+A_{n-1,n}x_n \\ D_nx_n\end{bmatrix}=\begin{bmatrix}0\\ \vdots\\ 0\\ 0\\  \vdots\\ 0\\ 0\end{bmatrix}.
$$
The last $n-k$ equalities give that $x_s=0$ for $s=k+1,..., n$ by the same reasoning seen in Case 1, which proves the claim. Therefore, $\alpha(T_n^d(A))=\sum\limits_{s=1}^k\alpha(D_s)<\infty$.

Secondly, we show that $\mathcal{R}(T_n^d(A))$ is closed. It is not hard to see that \\$\mathcal{R}(T_n^d(A))=\mathcal{R}(D_1)\oplus\cdots\oplus\mathcal{R}(D_{k-1})\oplus\mathcal{R}(D_k)\oplus\mathcal{R}(J_k)\oplus\cdots\oplus\mathcal{R}(D_{n-1})\oplus\mathcal{R}(J_{n-1})\oplus\mathcal{R}(D_n)$. Furthermore, due to left invertibility of $J_i$'s, there exist closed subspaces $U_i$ of $\mathcal{R}(D_i)'$ such that $\mathcal{R}(D_i)'=\mathcal{R}(J_i)\oplus U_i$, $k\leq i\leq n-1$. It means that $X_1\oplus X_2\oplus\cdots\oplus X_n=X_1\oplus\cdots\oplus X_{k-1}\oplus\mathcal{R}(D_{k})\oplus\mathcal{R}(J_k)\oplus U_k\oplus\cdots\oplus\mathcal{R}(D_{n-1})\oplus\mathcal{R}(J_{n-1})\oplus U_{n-1}\oplus X_n$. Comparing these equalities, we conclude that $\mathcal{R}(T_n^d(A))$ is closed.\\

\noindent\textbf{Case 3}: $0\in\Omega_n$\\

In this case $0\in\bigcap\limits_{s=1}^n\rho_{le}(D_s)$, which implies $\alpha(D_1),...,\alpha(D_n)<\infty$ ($D_1,...,D_n$ are left Fredholm). Choose trivially $A_{ij}=0$ for every $1\leq i<j\leq n$. Then $\mathcal{N}(T_n^d(A))\cong\mathcal{N}(D_1)\oplus\cdots\oplus\mathcal{N}(D_n)$ and $\mathcal{R}(T_n^d(A))=\mathcal{R}(D_1)\oplus\cdots\oplus\mathcal{R}(D_n)$, and since $X_1\oplus X_2\oplus\cdots\oplus X_n=\mathcal{R}(D_{1})\oplus\mathcal{R}(D_1)'\oplus\cdots\oplus\mathcal{R}(D_{n})\oplus\mathcal{R}(D_n')$, it follows that $\mathcal{R}(T_n^d(A))$ is closed and $\alpha(T_n^d(A))=\sum\limits_{s=1}^n\alpha(D_s)<\infty$. $\square$
\begin{Posledica}
Let $D_1\in\mathcal{B}(X_1),\ D_2\in\mathcal{B}(X_2),...,D_n\in\mathcal{B}(X_n)$ with all $D_i,\ i=1,...,n $ being regular. Additionally, suppose that we have 
$$
\mathcal{N}(D_k)\preceq X_j/\mathcal{R}(D_j)\quad for\  all\  n\geq k>j\geq 1.
$$
Then
$$
\bigcap\limits_{A\in\mathcal{B}_n}\sigma_{le}(T_n^d(A))=\sigma_{le}(D_1)\cup\Big(\bigcup\limits_{k=2}^n\Delta_k\Big),
$$
where
$$
\Delta_k:=\Big\lbrace\lambda\in\sigma_{le}(D_k):\  \sum\limits_{s=1}^{k-1}\beta(D_s-\lambda)<\infty\Big\rbrace,\ 2\leq k\leq n.
$$
\end{Posledica}
The following is a dual result of Theorem \ref{LEVIFREDHOLM}.
\begin{Teorema}\label{DESNIFREDHOLM}
Let $D_1\in\mathcal{B}(X_1),\ D_2\in\mathcal{B}(X_2),...,D_n\in\mathcal{B}(X_n)$ with all $D_i,\ i=1,...,n $ satisfying the complements condition. Additionally, suppose that we have 
\begin{equation}\label{DRUGA}
X_k/\mathcal{R}(D_k)\preceq\mathcal{N}(D_j)\quad for\  all\  1\leq k<j\leq n.
\end{equation}
Then
$$
\bigcap\limits_{A\in\mathcal{B}_n}\sigma_{re}(T_n^d(A))=\sigma_{re}(D_n)\cup\Big(\bigcup\limits_{k=1}^{n-1}\Delta_k\Big),
$$
where
$$
\Delta_k:=\Big\lbrace\lambda\in\sigma_{re}(D_k):\  \sum\limits_{s=k+1}^{n}\alpha(D_s-\lambda)<\infty\Big\rbrace,\ 1\leq k\leq n-1.
$$
\end{Teorema}
\textbf{Proof: }From Lemma \ref{VEZA} it follows that $\sigma_{re}(T_n^d(A))=\sigma_{le}(T_n^d(A)')$. Now, the proof immediately follows from Theorem \ref{LEVIFREDHOLM}, and with regards to (\ref{DUALNA}). $\square$

\begin{Posledica}
Let $D_1\in\mathcal{B}(X_1),\ D_2\in\mathcal{B}(X_2),...,D_n\in\mathcal{B}(X_n)$ with all $D_i,\ i=1,...,n $ being regular. Additionally, suppose that we have 
$$
X_k/\mathcal{R}(D_k)\preceq \mathcal{N}(D_j)\quad for\  all\  1\leq k<j\leq n.
$$
Then
$$
\bigcap\limits_{A\in\mathcal{B}_n}\sigma_{re}(T_n^d(A))=\sigma_{re}(D_n)\cup\Big(\bigcup\limits_{k=1}^{n-1}\Delta_k\Big),
$$
where
$$
\Delta_k:=\Big\lbrace\lambda\in\sigma_{re}(D_k):\  \sum\limits_{s=k+1}^{n}\alpha(D_s-\lambda)<\infty\Big\rbrace,\ 1\leq k\leq n-1.
$$
\end{Posledica}

Now, we conclude with perturbation of essential spectrum. We put into use Definition \ref{JAKOUTAPANJE}.
\begin{Teorema}
\label{FREDHOLM}
Let $D_1\in\mathcal{B}(X_1),\ D_2\in\mathcal{B}(X_2),...,D_n\in\mathcal{B}(X_n)$ with all $D_i,\ i=1,...,n $ satisfying the complements condition. Additionally, suppose that we have 
\begin{equation}\label{TRECA}
\mathcal{N}(D_k)\preceq_s X_j/\mathcal{R}(D_j)\quad for\  all\  n\geq k>j\geq 1.
\end{equation}
Then
\begin{center}
$
\bigcap\limits_{A\in\mathcal{B}_n}\sigma_{e}(T_n^d(A))=\sigma_{le}(D_1)\bigcup\limits\sigma_{re}(D_n)\bigcup\limits\Big(\bigcup\limits_{k=2}^{n-1}\Delta_k\Big)\bigcup\limits\lbrace\lambda\in\mathds{C}:\ \beta(D_1-\lambda)=\infty,\ \sum\limits_{s=2}^n\alpha(D_s-\lambda)<\infty\rbrace\bigcup\limits\lbrace\lambda\in\mathds{C}:\ \alpha(D_n-\lambda)=\infty,\ \sum\limits_{s=1}^{n-1}\beta(D_s-\lambda)<\infty\rbrace,
$
\end{center}
where
\begin{center}
$
\Delta_k:=\Big\lbrace\lambda\in\sigma_{le}(D_k):\  \sum\limits_{s=1}^{k-1}\beta(D_s-\lambda)<\infty\Big\rbrace\bigcup\limits\Big\lbrace\lambda\in\sigma_{re}(D_k):\  \newline\sum\limits_{s=k+1}^{n}\alpha(D_s-\lambda)<\infty\Big\rbrace,\ 2\leq k\leq n-1.
$
\end{center}
\end{Teorema}
\textbf{Proof: }Write $\Delta_1:=\sigma_{le}(D_1)\bigcup\limits\sigma_{re}(D_n)$ and $\Delta_n:=\lbrace\lambda\in\mathds{C}:\ \beta(D_1-\lambda)=\infty,\ \sum\limits_{s=2}^n\alpha(D_s-\lambda)<\infty\rbrace\bigcup\limits\lbrace\lambda\in\mathds{C}:\ \alpha(D_n-\lambda)=\infty,\ \sum\limits_{s=1}^{n-1}\beta(D_s-\lambda)<\infty\rbrace.$ By Theorems \ref{LEVIFREDHOLM} and \ref{DESNIFREDHOLM} we clearly have $\bigcup\limits_{k=1}^n\Delta_k\subseteq\bigcap\limits_{A\in\mathcal{B}_n}\sigma_e(T_n^d(A))$. Without loss of generality, we only prove the case when $\lambda=0$ in what follows. \\
For the opposite inclusion, it suffices to prove that $0\not\in\bigcup\limits_{k=1}^n\Delta_k$ implies $0\not\in\bigcap\limits_{A\in\mathcal{B}_n}\sigma_e(T_n^d(A))$, i. e. there exists an operator tuple $A\in\mathcal{B}_n$ such that $T_n^d(A)$ is Fredholm.\\
Let $0\not\in\bigcup\limits_{k=1}^n\Delta_k$. Then, direct calculations show that $0\in\bigcup\limits_{l=1}^{n-1}\bigcup\limits_{k=2}^n\Omega_{lk}\ (l<k)$ or $0\in\bigcap\limits_{s=1}^n\rho_e(D_s)$, where\\
$
\Omega_{1k}=\Big\lbrace\lambda\in\bigcap\limits_{s=k+1}^n\rho_e(D_s)\bigcap\limits\rho_{le}(D_1)\bigcap\limits\rho_{re}(D_k):\\ 
\beta(D_1-\lambda)=\alpha(D_k-\lambda)=\infty\Big\rbrace,\ 2\leq k\leq n-1;\\
\Omega_{1n}=\lbrace\lambda\in\rho_{le}(D_1)\bigcap\limits\rho_{re}(D_n):\ 
\beta(D_1-\lambda)=\alpha(D_n-\lambda)=\infty\rbrace,\\
\Omega_{lk}=\Big\lbrace\lambda\in\bigcap\limits_{s=1}^{l-1}\rho_e(D_s)\bigcap\limits\rho_{le}(D_l)\bigcap\limits\rho_{re}(D_k)\bigcap\limits\bigcap\limits_{s=k+1}^n\rho_e(D_s):\\
\beta(D_l-\lambda)=\alpha(D_k-\lambda)=\infty\Big\rbrace,\ 2\leq l<k\leq n-1;\\
\Omega_{ln}=\Big\lbrace\lambda\in\bigcap\limits_{s=1}^{l-1}\rho_e(D_s)\bigcap\limits\rho_{le}(D_l)\bigcap\limits\rho_{re}(D_n):\\ 
\beta(D_l-\lambda)=\alpha(D_n-\lambda)=\infty\Big\rbrace,\ 2\leq l\leq n-1.
$\\[3mm]
There are three cases to consider:\\[3mm]
\textbf{Case 1: }$0\in\Omega_{1k},\ 2\leq k\leq n-1$. In this case $\beta(D_1)=\alpha(D_k)=\infty$, $\alpha(D_1),\beta(D_k)<\infty$ and $D_{k+1},...,D_n$ are Fredholm.
We find $A\in\mathcal{B}_n$ such that $\alpha(T_n^d(A))<\infty$ and $\beta(T_n^d(A))<\infty$.  We choose $A=(A_{ij})_{1\leq i<j\leq n}$ so that $A_{ij}=0$ if $j-i\neq 1$ or $j\neq k$, that is we place all nonzero operators of tuple $A$ on the superdiagonal or in the k-th column. It remains to define $A_{ij}$ for $j-i=1$, $1\leq i<j\leq n$, and for $j=k$. First notice that $A_{i,i+1}:X_{i+1}\rightarrow X_i$. Since all of diagonal entries satisfy the complements condition, we know that $X_{i+1}=\mathcal{N}(D_{i+1})\oplus\mathcal{N}(D_{i+1})'$, $X_i=\mathcal{R}(D_i)'\oplus\mathcal{R}(D_i)$, and due to (\ref{TRECA}) we get $\mathcal{N}(D_{i+1})\preceq_s\mathcal{R}(D_i)'$, and so there is left invertible operator $J_{i}:\mathcal{N}(D_{i+1})\rightarrow\mathcal{R}(D_i)'$ with $\dim(\mathcal{R}(D_i)'/\mathcal{R}(J_i))<\infty$. We put $A_{i,i+1}=\begin{bmatrix}J_{i} & 0\\ 0 & 0\end{bmatrix}:\begin{bmatrix}\mathcal{N}(D_{i+1})\\ \mathcal{N}(D_{i+1})'\end{bmatrix}\rightarrow\begin{bmatrix}\mathcal{R}(D_i)'\\ \mathcal{R}(D_i)\end{bmatrix}$, and we implement this procedure for all $1\leq i\leq k-2$. If $j=k$, then the reasoning analogous to the previous one gives us left invertible $J_i':\mathcal{N}(D_k)\rightarrow\mathcal{R}(D_i)'$ with $\dim(\mathcal{R}(D_i)'/\mathcal{R}(J_i'))<\infty$, and so we put $A_{ik}=\begin{bmatrix}J_i' & 0\\ 0 & 0\\ \end{bmatrix}:\begin{bmatrix}\mathcal{N}(D_{k})\\ \mathcal{N}(D_{k})'\end{bmatrix}\rightarrow\begin{bmatrix}\mathcal{R}(D_i)'\\ \mathcal{R}(D_i)\end{bmatrix}$ for all $1\leq i\leq k-1$. Finally, we place $A_{ij}=0$ if $j>k.$

Now we have chosen our $A$, we show that $\mathcal{N}(T_n^d(A))\cong\mathcal{N}(D_1)\oplus\mathcal{N}(D_{k+1})\oplus\cdots\oplus\mathcal{N}(D_n)$. Let us put $T_n^d(A)x=0$, where $x=x_1+\cdots+x_n\in X_1\oplus\cdots\oplus X_n$. The previous equality is then equivalent to 
$$\begin{bmatrix}D_1x_1+A_{12}x_2+A_{1k}x_k\\ D_2x_2+A_{23}x_3+A_{2k}x_{k}\\ \vdots\\ D_{k-2}x_{k-2}+A_{k-2,k-1}x_{k-1}+A_{k-2,k}x_{k}\\ D_{k-1}x_{k-1}+A_{k-1,k}x_k\\ D_kx_k\\D_{k+1}x_{k+1}\\ \vdots\\ D_nx_n\end{bmatrix}=\begin{bmatrix}0\\ 0\\ \vdots\\ 0\\ 0\\0\\0 \\ \vdots \\0\end{bmatrix}.
$$ 
The k-th equality yields $x_k\in\mathcal{N}(D_k)$. Therefore, by the procedure already seen, the (k-1)-th equality gives $x_{k-1}=x_k=0$. Now, the last system of equations implies  
$$\begin{bmatrix}D_1x_1+A_{12}x_2\\ D_2x_2+A_{23}x_3\\ \vdots\\ D_{k-2}x_{k-2}\\ 0\end{bmatrix}=\begin{bmatrix}0\\ 0\\ \vdots\\ 0\\ 0\end{bmatrix}.
$$ 
It is not hard to see that continuing this procedure iteratively, we get $x_k=x_{k-1}=\cdots=x_{2}=0$, and this proves our statement. This means that $\alpha(T_n^d(A))<\infty.$

Now, we prove that $\beta(T_n^d(A))<\infty$. We know that  $\mathcal{R}(T_n^d(A))=\mathcal{R}(D_1)\oplus\mathcal{R}(J_{1})\oplus\mathcal{R}(J_1')\oplus\mathcal{R}(D_2)\oplus\mathcal{R}(J_{2})\oplus\mathcal{R}(J_2')\oplus\cdots\oplus\mathcal{R}(D_{k-2})\oplus\mathcal{R}(J_{k-2})\oplus\mathcal{R}(J_{k-2}')\oplus\mathcal{R}(D_{k-1})\oplus\mathcal{R}(J_{k-1}')\oplus\mathcal{R}(D_k)\oplus\mathcal{R}(D_{k+1})\oplus\cdots\oplus\mathcal{R}(D_n)$. Furthermore, we know that there exist closed finite dimensional subspaces $U_i,\ U_i'$ of $\mathcal{R}(D_i)'$ such that $\mathcal{R}(D_i)'=\mathcal{R}(J_i)\oplus U_i=\mathcal{R}(J_i')\oplus U_i'$, $1\leq i\leq k-1$. It means that $X_1\oplus X_2\oplus\cdots\oplus X_n=\mathcal{R}(D_1)\oplus\mathcal{R}(J_{1})\oplus U_1\oplus\mathcal{R}(J_{1}')\oplus U_1'\oplus\mathcal{R}(D_2)\oplus\mathcal{R}(J_{2})\oplus U_2\oplus\mathcal{R}(J_{2}')\oplus U_2'\oplus\cdots\oplus\mathcal{R}(D_{k-2})\oplus\mathcal{R}(J_{k-2})\oplus U_{k-2}\oplus\mathcal{R}(J_{k-2}')\oplus U_{k-2}'\oplus\mathcal{R}(D_{k-1})\oplus\mathcal{R}(J_{k-1}')\oplus U_{k-1}'\oplus X_k\oplus X_{k+1}\oplus\cdots\oplus X_n$. Comparing these equalities, we conclude that $\dim((X_1\oplus\cdots\oplus X_n)/\mathcal{R}(T_n^d(A)))=\dim(U_1)+\dim(U_1')+\dim(U_2)+\dim(U_2')+\cdots+\dim(U_{k-2})+\dim(U_{k-2}')+\dim(U_{k-1}')+\beta(D_k)+\beta(D_{k+1})+\cdots+\beta(D_n)<\infty$.\\[3mm]

\textbf{Case 2:} $0\in\Omega_{1n}$

 In this case $\alpha(D_1)<\infty,\ \beta(D_n)<\infty$, and $\beta(D_1)=\alpha(D_n)=\infty$. We choose $A\in\mathcal{B}_n$ analogously as in Case 1. We shall get $\mathcal{N}(T_n^d(A))\cong\mathcal{N}(D_1)$, and therefore $\alpha(T_n^d(A))<\infty$. We also get $\dim((X_1\oplus\cdots\oplus X_n)/\mathcal{R}(T_n^d(A)))=\dim(U_1)+\dim(U_1')+\dim(U_2)+\dim(U_2')+\cdots+\dim(U_{n-2})+\dim(U_{n-2}')+\dim(U_{n-1}')+\beta(D_n)<\infty$, with subspaces $U_i,\ U_i'$ as in Case 1.\\[3mm]

\textbf{Case 3}: $0\in\Omega_{lk}$, $2\leq l<k\leq n-1.$

In this case $\beta(D_l)=\alpha(D_k)=\infty$, $\alpha(D_l),\beta(D_k)<\infty$ and $D_{k+1},...,D_n$ are Fredholm, as well as $D_1,...,D_{l-1}$.
We find $A\in\mathcal{B}_n$ such that $\alpha(T_n^d(A))<\infty$ and $\beta(T_n^d(A))<\infty$.  We choose $A=(A_{ij})_{1\leq i<j\leq n}$ so that $A_{ij}=0$ if $j-i\neq 1$ or $j\neq k$, that is we place all nonzero operators of tuple $A$ on the superdiagonal or in the k-th column. It remains to define $A_{ij}$ for $j-i=1$, $1\leq i<j\leq n$, and for $j=k$. First notice that $A_{i,i+1}:X_{i+1}\rightarrow X_i$. Since all of diagonal entries satisfy the complements condition, we know that $X_{i+1}=\mathcal{N}(D_{i+1})\oplus\mathcal{N}(D_{i+1})'$, $X_i=\mathcal{R}(D_i)'\oplus\mathcal{R}(D_i)$, and due to (\ref{TRECA}) we get $\mathcal{N}(D_{i+1})\preceq_s\mathcal{R}(D_i)'$, and so there is left invertible operator $J_{i}:\mathcal{N}(D_{i+1})\rightarrow\mathcal{R}(D_i)'$ with $\dim(\mathcal{R}(D_i)'/\mathcal{R}(J_i))<\infty$. We put $A_{i,i+1}=\begin{bmatrix}J_{i} & 0\\ 0 & 0\end{bmatrix}:\begin{bmatrix}\mathcal{N}(D_{i+1})\\ \mathcal{N}(D_{i+1})'\end{bmatrix}\rightarrow\begin{bmatrix}\mathcal{R}(D_i)'\\ \mathcal{R}(D_i)\end{bmatrix}$, and we implement this procedure for all $l\leq i\leq k-2$. If $j=k$, then the reasoning analogous to the previous one gives us left invertible $J_i':\mathcal{N}(D_k)\rightarrow\mathcal{R}(D_i)'$ with $\dim(\mathcal{R}(D_i)'/\mathcal{R}(J_i'))<\infty$, and so we put $A_{ik}=\begin{bmatrix}J_i' & 0\\ 0 & 0\\ \end{bmatrix}:\begin{bmatrix}\mathcal{N}(D_{k})\\ \mathcal{N}(D_{k})'\end{bmatrix}\rightarrow\begin{bmatrix}\mathcal{R}(D_i)'\\ \mathcal{R}(D_i)\end{bmatrix}$ for all $l\leq i\leq k-1$. Finally, we place $A_{ij}=0$ if $j>k.$

Now we have chosen our $A$, we show that $\mathcal{N}(T_n^d(A))\cong\mathcal{N}(D_1)\oplus\cdots\oplus\mathcal{N}(D_{l})\oplus\mathcal{N}(D_{k+1})\oplus\cdots\oplus\mathcal{N}(D_n)$. Let us put $T_n^d(A)x=0$, where $x=x_1+\cdots+x_n\in X_1\oplus\cdots\oplus X_n$. The previous equality is then equivalent to 
$$\begin{bmatrix}D_1x_1+A_{1k}x_k\\ D_2x_2+A_{2k}x_{k}\\ \vdots\\ D_{l-1}x_{l-1}+A_{l-1,k}x_{k}\\ D_{l}x_{l}+A_{l,l+1}x_{l+1}+A_{l,k}x_k\\ \vdots\\  D_{k-2}x_{k-2}+A_{k-2,k-1}x_{k-1}+A_{k-1,k}x_k\\ D_{k-1}x_{k-1}+A_{k-1,k}x_k\\ D_kx_k\\D_{k+1}x_{k+1}\\ \vdots\\ D_nx_n\end{bmatrix}=\begin{bmatrix}0\\ 0\\ \vdots\\ 0\\ 0\\ \vdots\\ 0\\ 0\\0\\0 \\ \vdots \\0\end{bmatrix}.
$$ 
By the same reasoning as seen in Case 1 we get $x_k=0$. Now, the last system of equations is equivalent to  
$$\begin{bmatrix}D_1x_1\\ D_2x_2\\ \vdots\\ D_{l-1}x_{l-1}\\ D_{l}x_{l}+A_{l,l+1}x_{l+1}\\ \vdots\\  D_{k-2}x_{k-2}+A_{k-2,k-1}x_{k-1}\\ D_{k-1}x_{k-1}\\ 0\\D_{k+1}x_{k+1}\\ \vdots\\ D_nx_n\end{bmatrix}=\begin{bmatrix}0\\ 0\\ \vdots\\ 0\\ 0\\ \vdots\\ 0\\ 0\\0\\0 \\ \vdots \\0\end{bmatrix}.
$$ 
It is not hard to see that continuing this reasoning iteratively, we get $x_k=x_{k-1}=\cdots=x_{l+1}=0$, and this proves our statement. This means that $\alpha(T_n^d(A))<\infty.$

Now we prove that $\beta(T_n^d(A))<\infty$. We know that  $\mathcal{R}(T_n^d(A))=\mathcal{R}(D_1)\oplus\mathcal{R}(J_1')\oplus\mathcal{R}(D_2)\oplus\mathcal{R}(J_2')\oplus\cdots\oplus\mathcal{R}(D_{l-1})\oplus\mathcal{R}(J_{l-1}')\oplus\mathcal{R}(D_{l})\oplus\mathcal{R}(J_{l})\oplus\mathcal{R}(J_{l}')\oplus\cdots\oplus\mathcal{R}(D_{k-2})\oplus\mathcal{R}(J_{k-2})\oplus\mathcal{R}(J_{k-2}')\oplus\mathcal{R}(D_{k-1})\oplus\mathcal{R}(J_{k-1}')\oplus\mathcal{R}(D_k)\oplus\mathcal{R}(D_{k+1})\oplus\cdots\oplus\mathcal{R}(D_n)$. Furthermore, we know that there exist closed finite dimensional subspaces $U_i,\ U_i'$ of $\mathcal{R}(D_i)'$ such that $\mathcal{R}(D_i)'=\mathcal{R}(J_i)\oplus U_i=\mathcal{R}(J_i')\oplus U_i'$ with $\dim U_i,\dim U_i'<\infty$, $l\leq i\leq k-1$. It means that $X_1\oplus X_2\oplus\cdots\oplus X_n=\mathcal{R}(D_1)\oplus\mathcal{R}(J_{1}')\oplus U_1'\oplus\mathcal{R}(D_2)\oplus\mathcal{R}(J_{2}')\oplus U_2'\oplus\cdots\oplus\mathcal{R}(D_{l-1})\oplus\mathcal{R}(J_{l-1}')\oplus U_{l-1}'\oplus\mathcal{R}(D_{l})\oplus\mathcal{R}(J_{l})\oplus U_{l}\oplus\mathcal{R}(J_{l}')\oplus U_{l'}\oplus\cdots\oplus\mathcal{R}(D_{k-2})\oplus\mathcal{R}(J_{k-2})\oplus U_{k-2}\oplus\mathcal{R}(J_{k-2}')\oplus U_{k-2}'\oplus\mathcal{R}(D_{k-1})\oplus\mathcal{R}(J_{k-1}')\oplus U_{k-1}'\oplus X_k\oplus X_{k+1}\oplus\cdots\oplus X_n$. Comparing these equalities we conclude that $\dim((X_1\oplus\cdots\oplus X_n)/\mathcal{R}(T_n^d(A)))=\dim(U_1')+\dim(U_2')+\cdots+\dim(U_{l-1}')+\dim(U_{l})+\dim(U_{l}')+\cdots+\dim(U_{k-2})+\dim(U_{k-2}')+\dim(U_{k-1}')+\beta(D_k)+\beta(D_{k+1})+\cdots+\beta(D_n)<\infty$.\\[3mm]

\textbf{Case 4}: $0\in\Omega_{ln}$

In this case $\beta(D_l)=\alpha(D_n)=\infty$, $\alpha(D_l),\beta(D_n)<\infty$ and $D_1,...,D_{l-1}$ are Fredholm. We choose $A\in\mathcal{B}_n$ analogously as in Case 3. Namely,  we choose $A=(A_{ij})_{1\leq i<j\leq n}$ so that $A_{ij}=0$ if $j-i\neq 1$ or $j\neq n$, that is we place all nonzero ope\-ra\-tors of tuple $A$ on the superdiagonal or in the n-th column. Following the same method as in the previous case, we shall get $\mathcal{N}(T_n^d(A))\cong\mathcal{N}(D_1)\oplus\cdots\oplus\mathcal{N}(D_l)$, and therefore $\alpha(T_n^d(A))<\infty$. We also get $\dim((X_1\oplus\cdots\oplus X_n)/\mathcal{R}(T_n^d(A)))=\dim(U_1)+\dim(U_1')+\dim(U_2)+\dim(U_2')+\cdots+\dim(U_{n-2})+\dim(U_{n-2}')+\dim(U_{n-1}')+\beta(D_n)<\infty$, with subspaces $U_i,\ U_i'$ as in Case 3.\\[3mm]

\textbf{Case 5}: $ 0\in\bigcap\limits_{s=1}^n\rho_e(D_s)$

Then we know $\alpha(D_1),...,\alpha(D_n)<\infty$, $\beta(D_1),...,\beta(D_n)<\infty$. Choose trivially $A_{ij}=0$ for all $1\leq i<j\leq n$ to get $\alpha(T_n^d(A))=\sum\limits_{s=1}^n\alpha(D_s)<\infty,\ \beta(T_n^d(A))=\sum\limits_{s=1}^n\beta(D_s)<\infty$. The proof is completed. $\square$
\begin{Posledica}
Let $D_1\in\mathcal{B}(X_1),\ D_2\in\mathcal{B}(X_2),...,D_n\in\mathcal{B}(X_n)$ with all $D_i,\ i=1,...,n $ being regular. Additionally, suppose that we have 
$$
\mathcal{N}(D_k)\preceq_s X_j/\mathcal{R}(D_j)\quad for\  all\  n\geq k>j\geq 1.
$$
Then
\begin{center}
$
\bigcap\limits_{A\in\mathcal{B}_n}\sigma_{e}(T_n^d(A))=\sigma_{le}(D_1)\bigcup\limits\sigma_{re}(D_n)\bigcup\limits\Big(\bigcup\limits_{k=2}^{n-1}\Delta_k\Big)\bigcup\limits\lbrace\lambda\in\mathds{C}:\ \beta(D_1-\lambda)=\infty,\ \sum\limits_{s=2}^n\alpha(D_s-\lambda)<\infty\rbrace\bigcup\limits\lbrace\lambda\in\mathds{C}:\ \alpha(D_n-\lambda)=\infty,\ \sum\limits_{s=1}^{n-1}\beta(D_s-\lambda)<\infty\rbrace,
$
\end{center}
where
\begin{center}
$
\Delta_k:=\Big\lbrace\lambda\in\sigma_{le}(D_k):\  \sum\limits_{s=1}^{k-1}\beta(D_s-\lambda)<\infty\Big\rbrace\bigcup\limits\Big\lbrace\lambda\in\sigma_{re}(D_k):\  \newline\sum\limits_{s=k+1}^{n}\alpha(D_s-\lambda)<\infty\Big\rbrace,\ 2\leq k\leq n-1.
$
\end{center}
\end{Posledica}
Having in mind the statement of Remark \ref{SLICNISU}, we conclude that the following \\Theorem also holds. Its proof is similar to the proof of Theorem \ref{FREDHOLM}, but technically much more demanding.
\begin{Teorema}
Let $D_1\in\mathcal{B}(X_1),\ D_2\in\mathcal{B}(X_2),...,D_n\in\mathcal{B}(X_n)$ with all $D_i,\ i=1,...,n $ satisfying the complements condition. Additionally, suppose that we have 
$$
\mathcal{N}(D_k)\cong X_j/\mathcal{R}(D_j)\ (u.f.d.s)\quad for\  all\  n\geq k>j\geq 1.
$$
Then
\begin{center}
$
\bigcap\limits_{A\in\mathcal{B}_n}\sigma_{e}(T_n^d(A))=\sigma_{le}(D_1)\bigcup\limits\sigma_{re}(D_n)\bigcup\limits\Big(\bigcup\limits_{k=2}^{n-1}\Delta_k\Big)\bigcup\limits\lbrace\lambda\in\mathds{C}:\ \beta(D_1-\lambda)=\infty,\ \sum\limits_{s=2}^n\alpha(D_s-\lambda)<\infty\rbrace\bigcup\limits\lbrace\lambda\in\mathds{C}:\ \alpha(D_n-\lambda)=\infty,\ \sum\limits_{s=1}^{n-1}\beta(D_s-\lambda)<\infty\rbrace,
$
\end{center}
where
\begin{center}
$
\Delta_k:=\Big\lbrace\lambda\in\sigma_{le}(D_k):\  \sum\limits_{s=1}^{k-1}\beta(D_s-\lambda)<\infty\Big\rbrace\bigcup\limits\Big\lbrace\lambda\in\sigma_{re}(D_k):\  \sum\limits_{s=k+1}^{n}\alpha(D_s-\lambda)<\infty\Big\rbrace,\ 2\leq k\leq n-1.
$
\end{center}
\end{Teorema}
\begin{Posledica}
Let $D_1\in\mathcal{B}(X_1),\ D_2\in\mathcal{B}(X_2),...,D_n\in\mathcal{B}(X_n)$ with all $D_i,\ i=1,...,n $ being regular. Additionally, suppose that we have 
$$
\mathcal{N}(D_k)\cong X_j/\mathcal{R}(D_j)\ (u.f.d.s)\quad for\  all\  n\geq k>j\geq 1.
$$
Then
\begin{center}
$
\bigcap\limits_{A\in\mathcal{B}_n}\sigma_{e}(T_n^d(A))=\sigma_{le}(D_1)\bigcup\limits\sigma_{re}(D_n)\bigcup\limits\Big(\bigcup\limits_{k=2}^{n-1}\Delta_k\Big)\bigcup\limits\lbrace\lambda\in\mathds{C}:\ \beta(D_1-\lambda)=\infty,\ \sum\limits_{s=2}^n\alpha(D_s-\lambda)<\infty\rbrace\bigcup\limits\lbrace\lambda\in\mathds{C}:\ \alpha(D_n-\lambda)=\infty,\ \sum\limits_{s=1}^{n-1}\beta(D_s-\lambda)<\infty\rbrace,
$
\end{center}
where
\begin{center}
$
\Delta_k:=\Big\lbrace\lambda\in\sigma_{le}(D_k):\  \sum\limits_{s=1}^{k-1}\beta(D_s-\lambda)<\infty\Big\rbrace\bigcup\limits\Big\lbrace\lambda\in\sigma_{re}(D_k):\  \sum\limits_{s=k+1}^{n}\alpha(D_s-\lambda)<\infty\Big\rbrace,\ 2\leq k\leq n-1.
$
\end{center}
\end{Posledica}

\subsection{When does $\sigma_*(T_n^d(A))=\bigcup\limits_{s=1}^n\sigma_*(D_s)$ hold for arbitrary $A\in\mathcal{B}_n$, $\sigma_*\in\lbrace \sigma_{le},\sigma_{re}, \sigma_e\rbrace$?}
In this subsection we provide an answer to Question 3 stated in Section 2. This question has been solved for the case of infinite dimensional separable Hilbert spaces by Wu and Huang in \cite{WU}. It turns out that their results stated in Section 3 of \cite{WU} hold, word by word, in arbitrary Banach spaces, with the only distinction that in every claim we have to assume validity of the complements condition for all $D_i$'s. We invite readers to copy the statements of Theorems 4, 5, 6 from \cite{WU}t and their Corollaries, by including the complements condition assumption and substituting Hilbert spaces $\mathcal{H}_1$,..., $\mathcal{H}_n$ by Banach spaces $X_1,...,X_n$. As of special interest, we recommend studying Corollaries 2, 6 and 10 in a connection with an interesting notion of he single valued extension property of Finch (\cite{FINC}). However, for the convenience of a reader, we decided to provide a couple of mentioned statements for illustration, and so we give analogues of Corollaries 4, 8 and 12 from \cite{WU}. Their proofs can be found in \cite{WU}, Section 3.
\begin{Posledica}
Let $D_1\in\mathcal{B}(X_1),\ D_2\in\mathcal{B}(X_2),...,D_n\in\mathcal{B}(X_n)$ with all $D_i,\ i=1,...,n $ satisfying the complements condition. Then
$$
\sigma_{le}(T_n^d(A))=\bigcup\limits_{k=1}^n\sigma_{le}(D_k)
$$
holds for every $A\in\mathcal{B}_n$ if and only if
$$
\begin{aligned}
\bigcup\limits_{k=2}^n\lbrace\lambda\in\sigma_{le}(D_k):\ \sum\limits_{s=1}^{k-1}\beta(D_s-\lambda)=\infty\rbrace\subseteq
\sigma_{le}(D_1)\cup\\ \bigcup\limits_{k=2}^{n-1}\lbrace\lambda\in\sigma_{le}(D_k):\ \sum\limits_{s=1}^{k-1}\beta(D_s-\lambda)<\infty\rbrace.
\end{aligned}
$$
\end{Posledica}
\begin{Posledica}
Let $D_1\in\mathcal{B}(X_1),\ D_2\in\mathcal{B}(X_2),...,D_n\in\mathcal{B}(X_n)$ with all $D_i,\ i=1,...,n $ satisfying the complements condition. Then
$$
\sigma_{re}(T_n^d(A))=\bigcup\limits_{k=1}^n\sigma_{re}(D_k)
$$
holds for every $A\in\mathcal{B}_n$ if and only if
$$
\begin{aligned}
\bigcup\limits_{k=1}^{n-1}\lbrace\lambda\in\sigma_{re}(D_k):\ \sum\limits_{s=k+1}^{n}\alpha(D_s-\lambda)=\infty\rbrace\subseteq
\sigma_{re}(D_n)\cup\\ \bigcup\limits_{k=2}^{n-1}\lbrace\lambda\in\sigma_{re}(D_k):\ \sum\limits_{s=k+1}^{n}\alpha(D_s-\lambda)<\infty\rbrace.
\end{aligned}
$$
\end{Posledica}
\begin{Posledica}
Let $D_1\in\mathcal{B}(X_1),\ D_2\in\mathcal{B}(X_2),...,D_n\in\mathcal{B}(X_n)$ with all $D_i,\ i=1,...,n $ satisfying the complements condition. Then
$$
\sigma_{e}(T_n^d(A))=\bigcup\limits_{k=1}^n\sigma_{e}(D_k)
$$
holds for every $A\in\mathcal{B}_n$ if and only if
$$
\begin{aligned}
\bigcup\limits_{k=1}^{n-1}\lbrace\lambda\in\mathds{C}:\ \beta(D_k-\lambda)=\sum\limits_{s=k+1}^{n}\alpha(D_s-\lambda)=\infty\rbrace\subseteq\\
\sigma_{le}(D_1)\cup\sigma_{re}(D_n)\cup\Big(\bigcup\limits_{k=2}^{n-1}\Delta_k\Big),
\end{aligned}
$$
where $\Delta_k$ is defined in Theorem \ref{FREDHOLM}.
\end{Posledica}

\section{Perturbations of left (right) the Weyl spectrum of $T_n^d(A)$}

We start with a result which characterizes perturbations of left the Weyl spectrum of $T_n^d(A)$.  We need to use stronger form of Definition \ref{UTAPANJE}, namely Definition \ref{ESENCIJALNOUTAPANJE}.
\begin{Teorema}\label{LEVIVEJL}
Let $D_1\in\mathcal{B}(X_1),\ D_2\in\mathcal{B}(X_2),...,D_n\in\mathcal{B}(X_n)$ with all $D_i,\ i=1,...,n $ satisfying the complements condition. Additionally, suppose that we have 
\begin{equation}\label{CETVRTA}
\mathcal{N}(D_k)\prec X_j/\mathcal{R}(D_j)\quad for\  all\  n\geq k>j\geq 1.
\end{equation}
Then
$$
\bigcap\limits_{A\in\mathcal{B}_n}\sigma_{lw}(T_n^d(A))=\sigma_{le}(D_1)\cup\Big(\bigcup\limits_{k=2}^{n+1}\Delta_k\Big),
$$
where
$$
\Delta_k:=\Big\lbrace\lambda\in\sigma_{le}(D_k):\  \sum\limits_{s=1}^{k-1}\beta(D_s-\lambda)<\infty\Big\rbrace,\ 2\leq k\leq n,
$$
$$
\Delta_{n+1}:=\lbrace\lambda\in\mathds{C}:\ \sum\limits_{s=1}^n\beta(D_s-\lambda)<\sum\limits_{s=1}^n\alpha(D_s-\lambda)\rbrace.
$$
\end{Teorema}
\textbf{Proof. }Inclusion $(\supseteq)$ is already proved in \cite{WU2}. For convenience of a reader we provide the full proof, with a slight difference in that now we have different space decompositions.

Write $\Delta_1:=\sigma_{le}(D_1)$,
$
\Delta_{k1}=\lbrace\lambda\in\mathds{C}: \mathcal{R}(D_k-\lambda)\  is\  not\  closed\ and\\ \ \sum\limits_{s=1}^{k-1}\beta(D_s-\lambda)<\infty\rbrace
$
and
$
\Delta_{k2}=\lbrace\lambda\in\mathds{C}: \mathcal{R}(D_k-\lambda)\  is\  closed\  and\\  \sum\limits_{s=1}^{k-1}\beta(D_s-\lambda)<\infty,\ \alpha(D_k-\lambda)=\infty\rbrace$
for $2\leq k\leq n$. \\Lemma \ref{POMOCNALEMA} implies $\Delta_1\subseteq\bigcap\limits_{A\in\mathcal{B}_n}\sigma_{lw}(T_n^d(A))$ . Let $\lambda\in\bigcup\limits_{k=2}^n\Delta_k$.

Let $\lambda\in\Delta_{k1}$, that is $\mathcal{R}(D_k-\lambda)$  is not closed and $\sum\limits_{s=1}^{k-1}\beta(D_s-\lambda)<\infty$, $k\in\lbrace2,...,n\rbrace$. Then for each $A\in\mathcal{B}_k$, operator matrix $T_k^d(A)$ as an operator from $\mathcal{N}(D_1)'\oplus\mathcal{N}(D_1)\oplus\mathcal{N}(D_2)'\oplus\mathcal{N}(D_2)\oplus\mathcal{N}(D_3)'\oplus\mathcal{N}(D_3)\oplus\cdots\oplus\mathcal{N}(D_k)'\oplus\mathcal{N}(D_k)$ into $\mathcal{R}(D_1)\oplus\mathcal{R}(D_1)'\oplus\mathcal{R}(D_2)\oplus\mathcal{R}(D_2)'\oplus\cdots\oplus\mathcal{R}(D_{k-1})\oplus\mathcal{R}(D_{k-1}')\oplus\mathcal{R}(D_k)\oplus\mathcal{R}(D_k)'$ admits the following block representation
\begin{equation}\label{MATRICA'}
T_k^d(A)=\scalemath{0.95}
{\begin{bmatrix}
    D_1^{(1)} & 0 & A_{12}^{(1)} & A_{12}^{(2)} & ... & A_{1,k-1}^{(1)} & A_{1,k-1}^{(2)} & A_{1k}^{(1)} & A_{1k}^{(2)}\\
    0 & 0 & A_{12}^{(3)} & A_{12}^{(4)} & ... & A_{1,k-1}^{(3)} & A_{1,k-1}^{(4)} & A_{1k}^{(3)} & A_{1k}^{(4)}\\
    0 & 0 & D_2^{(1)} & 0 & ... & A_{2,k-1}^{(1)} & A_{2,k-1}^{(2)} & A_{2k}^{(1)} & A_{2k}^{(2)}\\
    0 & 0 & 0 & 0 & ... & A_{2,k-1}^{(3)} & A_{2,k-1}^{(4)} & A_{2k}^{(3)} & A_{2k}^{(4)}\\
    \vdots & \vdots & \vdots & \vdots & \ddots & \vdots & \vdots & \vdots & \vdots\\
    0 & 0 & 0 & 0 & ... & D_{k-1}^{(1)} & 0 & A_{k-1,k}^{(1)} & A_{k-1,k}^{(2)}\\
    0 & 0 & 0 & 0 & ... & 0 & 0 & A_{k-1,k}^{(3)} & A_{k-1,k}^{(4)}\\
    0 & 0 & 0 & 0 & ... & 0 & 0 & D_k^{(1)} & 0\\
    0 & 0 & 0 & 0 & ... & 0 & 0 & 0 & 0\\
\end{bmatrix}}
\end{equation}
Assume that $\mathcal{R}(T_k^d(A)-\lambda)$ is closed. Then, according to Lemma \ref{PRVALEMA}, range of the following operator matrix
$$
\begin{bmatrix}
    D_1^{(1)}-\lambda & 0 & A_{12}^{(1)} & A_{12}^{(2)} & ... & A_{1,k-1}^{(1)} & A_{1,k-1}^{(2)} & A_{1k}^{(1)} & A_{1k}^{(2)}\\
    0 & 0 & D_2^{(1)}-\lambda & 0 & ... & A_{2,k-1}^{(1)} & A_{2,k-1}^{(2)} & A_{2k}^{(1)} & A_{2k}^{(2)}\\
    \vdots & \vdots & \vdots & \vdots & \ddots & \vdots & \vdots & \vdots & \vdots\\
    0 & 0 & 0 & 0 & ... & D_{k-1}^{(1)}-\lambda & 0 & A_{k-1,k}^{(1)} & A_{k-1,k}^{(2)}\\
    0 & 0 & 0 & 0 & ... & 0 & 0 & D_k^{(1)}-\lambda & 0\\
\end{bmatrix}
$$
is closed due to $\sum\limits_{s=1}^{k-1}\beta(D_s-\lambda)<\infty$. Since $\overline{\mathcal{R}(D_s^{(1)}-\lambda)}=\overline{\mathcal{R}(D_s-\lambda)}$, $s=1,2,...,k-1$
, then
$$
\scalemath{0.95}{\mathcal{R}\left(\begin{bmatrix}D_1^{(1)}-\lambda & A_{12}^{(1)} & A_{12}^{(2)} & ... & A_{1,k-1}^{(1)} & A_{1,k-1}^{(2)}\\
    0 & D_2^{(1)}-\lambda & 0 & ... & A_{2,k-1}^{(1)} & A_{2,k-1}^{(2)}\\
    \vdots & \vdots & \vdots & \ddots & \vdots & \vdots\\
    0 & 0 & 0 & ... & D_{k-1}^{(1)}-\lambda & 0\\\end{bmatrix}\right)}
$$
is dense in $\overline{\mathcal{R}(D_1-\lambda)}\oplus\overline{\mathcal{R}(D_2-\lambda)}\oplus\cdots\oplus\overline{\mathcal{R}(D_{k-1}-\lambda)}$ from Lemma \ref{VELIKALEMA}. By Lemma \ref{DRUGALEMA}, it follows that $\mathcal{R}(D_k^{(1)}-\lambda)$ is closed, which is contradiction to the fact that $\mathcal{R}(D_k-\lambda)$ is not closed. This implies that $\mathcal{R}(T_k^d(A)-\lambda)$ is not closed for every $A\in\mathcal{B}_k$, and hence $\lambda\in\bigcap\limits_{A\in\mathcal{B}_k}\sigma_{le}(T_k^d(A))$. This proves $\bigcup\limits_{k=2}^n\Delta_{k1}\subseteq\bigcap\limits_{A\in\mathcal{B}_n}\sigma_{lw}(T_n^d(A))$.

Let $\lambda\in\Delta_{k2}$, $k\in\lbrace 2,...,n\rbrace.$ Without loss of generality, put $\lambda=0$. We only consider $0\in\Delta_{k2}\setminus(\bigcup\limits_{s=2}^k\Delta_{s1}\bigcup\limits\Delta_1)$. Therefore, $0\in\rho_{le}(D_1)$, $\mathcal{R}(D_s)$ is closed for $s=2,...,k$, $\sum\limits_{s=1}^{k-1}\beta(D_s)<\infty$ and $\alpha(D_k)=\infty$. Then for each $A\in\mathcal{B}_k$, operator matrix $T_k^d(A)$ admits block representation as in (\ref{MATRICA'}). Evidently, $D_1^{(1)},\ D_2^{(1)},...,D_k^{(1)}$ as in (\ref{MATRICA'}) are all invertible. Thus, there exist invertible operator matrices $U$ and $V$ such that 
\begin{equation}\label{MATRICA''}
UT_k^d(A)V=\scalemath{0.85}{\begin{bmatrix}
D_1^{(1)} & 0 & 0 & 0 & ... & 0 & 0 & 0 & 0\\
0 & 0 & 0 & A_{12}^{(4)} & ... & 0 & A_{1,k-1}^{(4)} & 0 & A_{1k}^{(4)}\\
0 & 0 & D_2^{(1)} & 0 & ... & 0 & 0 & 0 & 0\\
0 & 0 & 0 & 0 & ... & 0 & A_{2,k-1}^{(4)} & 0 & A_{2k}^{(4)}\\
0 & 0 & 0 & 0 & ... & 0 & 0 & 0 & 0\\
0 & 0 & 0 & 0 & ... & 0 & A_{3,k-1}^{(4)} & 0 & A_{3k}^{(4)}\\
\vdots & \vdots & \vdots & \vdots & \ddots & \vdots & \vdots & \vdots & \vdots\\
0 & 0 & 0 & 0 & ... & D_{k-1}^{(4)} & 0 & 0 & 0\\
0 & 0 & 0 & 0 & ... & 0 & 0 & 0 & A_{k-1,k}^{(4)}\\
0 & 0 & 0 & 0 & ... & 0 & 0 & D_k^{(1)} & 0\\
0 & 0 & 0 & 0 & ... & 0 & 0 & 0 & 0\\
\end{bmatrix}
}
\end{equation}
Note that $A_{ij}^{(4)}$ in (\ref{MATRICA''}) are not the original ones in (\ref{MATRICA'}) in general, but we still use them for the convenience. Combine by (\ref{MATRICA''}), $T_k^d(A)$ is left Weyl if and only if
\begin{equation}\label{MATRICA'''}
\begin{bmatrix}
0 & A_{12}^{(4)} & A_{13}^{(4)} & A_{14}^{(4)} & ... & A_{1k}^{(4)}\\
0 & 0          & A_{23}^{(4)} & A_{24}^{(4)} & ... & A_{2k}^{(4)}\\
0 & 0           & 0         & A_{34}^{(4)} & ... & A_{3k}^{(4)}\\
\vdots & \vdots   &   \vdots & \vdots & \ddots & \vdots\\
0 & 0          &   0        &    0      & ... & A_{k-1,k}^{(4)}\\
0 & 0 & 0 & 0 & ... & 0\\
\end{bmatrix}
:
\begin{bmatrix}
\mathcal{N}(D_1)\\
\mathcal{N}(D_2)\\
\mathcal{N}(D_3)\\
\vdots\\
\mathcal{N}(D_{k-1})\\
\mathcal{N}(D_k)\\
\end{bmatrix}
\rightarrow
\begin{bmatrix}
\mathcal{R}(D_1)'\\
\mathcal{R}(D_2)'\\
\mathcal{R}(D_3)'\\
\vdots\\
\mathcal{R}(D_{k-1})'\\
\mathcal{R}(D_k)'\\
\end{bmatrix}
\end{equation}
is left Weyl. Since $\sum\limits_{s=1}^{k-1}\beta(D_s)<\infty$ and $\alpha(D_k)=\infty$, it follows that 
$$
\alpha\left(\begin{bmatrix}
A_{1k}^{(4)}\\
A_{2k}^{(4)}\\
A_{3k}^{(4)}\\
\vdots\\
A_{k-1,k}^{(4)}
\end{bmatrix}\right)=\infty,
$$
and hence the operator defined in (\ref{MATRICA'''}) is not left Weyl for every $A\in\mathcal{B}_k$. This proves $\bigcup\limits_{k=2}^n\Delta_{k2}\subseteq\bigcap\limits_{A\in\mathcal{B}_n}\sigma_{lw}(T_n^d(A))$.\\

Let $0\in\Delta_{n+1}\setminus\bigcup\limits_{s=1}^n\Delta_s$. Then $0\in\bigcap\limits_{s=1}^n\rho_{le}(D_s)$ and $\sum\limits_{s=1}^n\beta(D_s)<\sum\limits_{s=1}^n\alpha(D_s)$. Thus for each $A\in\mathcal{B}_n$, operator matrix $T_n^d(A)$ admits block representation as (\ref{MATRICA'}), and we still use (\ref{MATRICA''}) and (\ref{MATRICA'''}), where $k=n$. Since $0\in\bigcap\limits_{s=1}^n\rho_{le}(D_s)$, it follows that $T_n^d(A)$ is left Weyl if and only if operator (\ref{MATRICA'''}) is left Weyl. From $\sum\limits_{s=1}^n\beta(D_s)<\sum\limits_{s=1}^n\alpha(D_s)$, we know (\ref{MATRICA'''}) is not left Weyl for every $A\in\mathcal{B}_n$. Therefore, $\Delta_{n+1}\subseteq\bigcap\limits_{A\in\mathcal{B}_n}\sigma_{lw}(T_n^d(A))$.\\[2mm]

For the opposite inclusion, it suffices to prove that $0\not\in\bigcup\limits_{k=1}^{n+1}\Delta_k$ implies $0\not\in\bigcap\limits_{A\in\mathcal{B}_n}\sigma_{lw}(T_n^d(A))$, i.e. there exists an operator tuple $A\in\mathcal{B}_n$ such that $T_n^d(A)$ is left Weyl. Let $0\not\in\bigcup\limits_{k=1}^{n+1}\Delta_k$. Then, direct calculation show that $0\in\bigcup\limits_{k=1}^n\Omega_k$, where\\
$
\Omega_1=\lbrace\lambda\in\rho_{le}(D_1):\ \beta(D_1-\lambda)=\infty\rbrace\\
\Omega_k=\lbrace\lambda\in\bigcap\limits_{s=1}^k\rho_{le}(D_s):\ \beta(D_k-\lambda)=\infty\rbrace,\ 2\leq k\leq n-1,\\
\Omega_n=\lbrace\lambda\in\bigcap\limits_{s=1}^n\rho_{le}(D_s):\ \sum\limits_{s=1}^n\alpha(D_s-\lambda)\leq\sum\limits_{s=1}^n\beta(D_s-\lambda)\rbrace.\\
$There are three cases to consider.\\[3mm]
\textbf{Case 1}: $0\in\Omega_1$\\

This means that $\alpha(D_1)<\infty$ ($D_1$ is left Fredholm), and $\beta(D_1)=\infty$. We find $A\in\mathcal{B}_n$ such that $\alpha(T_n^d(A))<\infty$, $\beta(T_n^d(A))=\infty$ and $\mathcal{R}(T_n^d(A))$ is closed.  We choose $A=(A_{ij})_{1\leq i<j\leq n}$ completely the same as in Case 1 in the proof of Theorem \ref{LEVIFREDHOLM}. We get $\mathcal{N}(T_n^d(A))\cong\mathcal{N}(D_1)$ and $\mathcal{R}(T_n^d(A))=\mathcal{R}(D_1)\oplus\mathcal{R}(J_{1})\oplus\mathcal{R}(D_2)\oplus\mathcal{R}(J_{2})\oplus\cdots\oplus\mathcal{R}(D_{n-1})\oplus\mathcal{R}(J_{n-1})\oplus\mathcal{R}(D_n)$, \\$X_1\oplus X_2\oplus\cdots\oplus X_n=\mathcal{R}(D_1)\oplus\mathcal{R}(J_{1})\oplus U_1\oplus\mathcal{R}(D_2)\oplus\mathcal{R}(J_{2})\oplus U_2\oplus\cdots\oplus\mathcal{R}(D_{n-1})\oplus\mathcal{R}(J_{n-1})\oplus U_{n-1}\oplus\mathcal{R}(D_n)$. Therefore, $\alpha(T_n^d(A))<\infty$ and $\beta(T_n^d(A))=\dim(U_1)+\dim(U_2)+\cdots+\dim(U_{n-1}) +\beta(D_n)=\infty$ due to (\ref{CETVRTA}). It follows that $T_n^d(A)$ is left Weyl.\\

\noindent\textbf{Case 2}: $0\in\Omega_k$, $2\leq k\leq n-1$\\

This means that $\alpha(D_1),...,\alpha(D_k)<\infty$ ($D_1,..., D_k$ are left Fredholm), and $\beta(D_k)=\infty$. We find $A\in\mathcal{B}_n$ such that $\alpha(T_n^d(A))<\infty$, $\beta(T_n^d(A))=\infty$ and $\mathcal{R}(T_n^d(A))$ is closed.  We choose $A=(A_{ij})_{1\leq i<j\leq n}$ completely the same as in Case 2 in the proof of Theorem \ref{LEVIFREDHOLM}. We get $\alpha(T_n^d(A))=\sum\limits_{s=1}^k\alpha(D_s)<\infty$ and $\beta(T_n^d(A))=\beta(D_1)+\cdots+\beta(D_{k-1})+\dim(U_k)+\cdots+\dim(U_{n-1})+\beta(D_n)=\infty$ due to (\ref{CETVRTA}). Therefore $T_n^d(A)$ is left Weyl.\\

\noindent\textbf{Case 3}: $0\in\Omega_n$\\

In this case $0\in\bigcap\limits_{s=1}^n\rho_{le}(D_s)$ and $\sum\limits_{s=1}^n\alpha(D_s)-\sum\limits_{s=1}^n\beta(D_s)\leq0$, which implies $\alpha(D_1),...,\alpha(D_n)<\infty$ ($D_1,...,D_n$ are left Fredholm). Choose trivially $A_{ij}=0$ for every $1\leq i<j\leq n$. Then $\alpha(T_n^d(A))=\sum\limits_{s=1}^n\alpha(D_s)<\infty$, $\beta(T_n^d(A))=\sum\limits_{s=1}^n\beta(D_s)$. Therefore, $\alpha(T_n^d(A))-\beta(T_n^d(A))\leq0$, hence $T_n^d(A)$ is left Weyl. $\square$
\begin{Posledica}
Let $D_1\in\mathcal{B}(X_1),\ D_2\in\mathcal{B}(X_2),...,D_n\in\mathcal{B}(X_n)$ with all $D_i,\ i=1,...,n $ being regular. 
Additionally, suppose that we have 
\begin{equation}\label{CETVRTA}
\mathcal{N}(D_k)\prec X_j/\mathcal{R}(D_j)\quad for\  all\  n\geq k>j\geq 1.
\end{equation}
Then
$$
\bigcap\limits_{A\in\mathcal{B}_n}\sigma_{lw}(T_n^d(A))=\sigma_{le}(D_1)\cup\Big(\bigcup\limits_{k=2}^{n+1}\Delta_k\Big),
$$
where
$$
\Delta_k:=\Big\lbrace\lambda\in\sigma_{le}(D_k):\  \sum\limits_{s=1}^{k-1}\beta(D_s-\lambda)<\infty\Big\rbrace,\ 2\leq k\leq n,
$$
$$
\Delta_{n+1}:=\lbrace\lambda\in\mathds{C}:\ \sum\limits_{s=1}^n\beta(D_s-\lambda)<\sum\limits_{s=1}^n\alpha(D_s-\lambda)\rbrace.
$$
\end{Posledica}
Inspecting the proof of Theorem \ref{LEVIVEJL}, one easily spots that we can replace condition (\ref{CETVRTA}) with what seems to be looser condition. And so, we state the following.
\begin{Teorema}\label{LEVIVEJL'}
Let $D_1\in\mathcal{B}(X_1),\ D_2\in\mathcal{B}(X_2),...,D_n\in\mathcal{B}(X_n)$ with all $D_i,\ i=1,...,n $ satisfying the complements condition. Additionally, suppose that we have 
\begin{center}
$\mathcal{N}(D_k)\preceq X_j/\mathcal{R}(D_j)\quad for\  all\  n\geq k>j\geq 1$,\\[1mm]
$\beta(D_n)=\infty.$
\end{center}
Then
$$
\bigcap\limits_{A\in\mathcal{B}_n}\sigma_{lw}(T_n^d(A))=\sigma_{le}(D_1)\cup\Big(\bigcup\limits_{k=2}^{n+1}\Delta_k\Big),
$$
where
$$
\Delta_k:=\Big\lbrace\lambda\in\sigma_{le}(D_k):\  \sum\limits_{s=1}^{k-1}\beta(D_s-\lambda)<\infty\Big\rbrace,\ 2\leq k\leq n,
$$
$$
\Delta_{n+1}:=\lbrace\lambda\in\mathds{C}:\ \sum\limits_{s=1}^n\beta(D_s-\lambda)<\sum\limits_{s=1}^n\alpha(D_s-\lambda)\rbrace.
$$
\end{Teorema}
\begin{Posledica}
Let $D_1\in\mathcal{B}(X_1),\ D_2\in\mathcal{B}(X_2),...,D_n\in\mathcal{B}(X_n)$ with all $D_i,\ i=1,...,n $ being regular. 
Additionally, suppose that we have 
\begin{center}
$\mathcal{N}(D_k)\preceq X_j/\mathcal{R}(D_j)\quad for\  all\  n\geq k>j\geq 1$,\\[1mm]
$\beta(D_n)=\infty.$
\end{center}
Then
$$
\bigcap\limits_{A\in\mathcal{B}_n}\sigma_{lw}(T_n^d(A))=\sigma_{le}(D_1)\cup\Big(\bigcup\limits_{k=2}^{n+1}\Delta_k\Big),
$$
where
$$
\Delta_k:=\Big\lbrace\lambda\in\sigma_{le}(D_k):\  \sum\limits_{s=1}^{k-1}\beta(D_s-\lambda)<\infty\Big\rbrace,\ 2\leq k\leq n,
$$
$$
\Delta_{n+1}:=\lbrace\lambda\in\mathds{C}:\ \sum\limits_{s=1}^n\beta(D_s-\lambda)<\sum\limits_{s=1}^n\alpha(D_s-\lambda)\rbrace.
$$
\end{Posledica}
The following is a dual result of Theorem \ref{LEVIVEJL}.
\begin{Teorema}\label{DESNIVEJL}
Let $D_1\in\mathcal{B}(X_1),\ D_2\in\mathcal{B}(X_2),...,D_n\in\mathcal{B}(X_n)$ with all $D_i,\ i=1,...,n $ satisfying the complements condition. Additionally, suppose that we have 
\begin{equation}\label{PETA}
X_k/\mathcal{R}(D_k)\prec\mathcal{N}(D_j)\quad for\  all\  1\leq k<j\leq n.
\end{equation}
Then
$$
\bigcap\limits_{A\in\mathcal{B}_n}\sigma_{rw}(T_n^d(A))=\sigma_{re}(D_n)\cup\Big(\bigcup\limits_{k=1}^{n-1}\Delta_k\Big)\cup\Delta_{n+1},
$$
where
$$
\Delta_k:=\Big\lbrace\lambda\in\sigma_{re}(D_k):\  \sum\limits_{s=k+1}^{n}\alpha(D_s-\lambda)<\infty\Big\rbrace,\ 1\leq k\leq n-1.
$$
$$
\Delta_{n+1}=\lbrace\lambda\in\mathds{C}:\ \sum\limits_{s=1}^n\alpha(D_s-\lambda)<\sum\limits_{s=1}^n\beta(D_s-\lambda)\rbrace.
$$
\end{Teorema}
\textbf{Proof: }From Lemma \ref{VEZA} it follows that $\sigma_{rw}(T_n^d(A))=\sigma_{lw}(T_n^d(A)')$. Now, proof immediately follows from Theorem \ref{LEVIVEJL}, and with regards to the form (\ref{DUALNA}). $\square$

\begin{Posledica}
Let $D_1\in\mathcal{B}(X_1),\ D_2\in\mathcal{B}(X_2),...,D_n\in\mathcal{B}(X_n)$ with all $D_i,\ i=1,...,n $ being regular. Additionally, suppose that we have 
$$
X_k/\mathcal{R}(D_k)\prec \mathcal{N}(D_j)\quad for\  all\  1\leq k<j\leq n.
$$
Then
$$
\bigcap\limits_{A\in\mathcal{B}_n}\sigma_{rw}(T_n^d(A))=\sigma_{re}(D_n)\bigcup\limits\Big(\bigcup\limits_{k=1}^{n-1}\Delta_k\Big),
$$
where
$$
\Delta_k:=\Big\lbrace\lambda\in\sigma_{re}(D_k):\  \sum\limits_{s=k+1}^{n}\alpha(D_s-\lambda)<\infty\Big\rbrace,\ 1\leq k\leq n-1.
$$
$$
\Delta_{n+1}=\lbrace\lambda\in\mathds{C}:\ \sum\limits_{s=1}^n\alpha(D_s-\lambda)<\sum\limits_{s=1}^n\beta(D_s-\lambda)\rbrace.
$$
\end{Posledica}
\begin{Teorema}\label{DESNIVEJL'}
Let $D_1\in\mathcal{B}(X_1),\ D_2\in\mathcal{B}(X_2),...,D_n\in\mathcal{B}(X_n)$ with all $D_i,\ i=1,...,n $ satisfying the complements condition. Additionally, suppose that we have 
\begin{center}
$X_k/\mathcal{R}(D_k)\preceq \mathcal{N}(D_j)\quad for\  all\  1\leq k<j\leq n$,\\[1mm]
$\alpha(D_1)=\infty.$
\end{center}
Then
$$
\bigcap\limits_{A\in\mathcal{B}_n}\sigma_{rw}(T_n^d(A))=\sigma_{re}(D_n)\cup\Big(\bigcup\limits_{k=1}^{n-1}\Delta_k\Big)\cup\Delta_{n+1},
$$
where
$$
\Delta_k:=\Big\lbrace\lambda\in\sigma_{re}(D_k):\  \sum\limits_{s=k+1}^{n}\alpha(D_s-\lambda)<\infty\Big\rbrace,\ 1\leq k\leq n-1.
$$
$$
\Delta_{n+1}=\lbrace\lambda\in\mathds{C}:\ \sum\limits_{s=1}^n\alpha(D_s-\lambda)<\sum\limits_{s=1}^n\beta(D_s-\lambda)\rbrace.
$$
\end{Teorema}
\begin{Posledica}
Let $D_1\in\mathcal{B}(X_1),\ D_2\in\mathcal{B}(X_2),...,D_n\in\mathcal{B}(X_n)$ with all $D_i,\ i=1,...,n $ being regular. Additionally, suppose that we have 
\begin{center}
$X_k/\mathcal{R}(D_k)\preceq \mathcal{N}(D_j)\quad for\  all\  1\leq k<j\leq n$,\\[1mm]
$\alpha(D_1)=\infty.$
\end{center}
Then
$$
\bigcap\limits_{A\in\mathcal{B}_n}\sigma_{rw}(T_n^d(A))=\sigma_{re}(D_n)\cup\Big(\bigcup\limits_{k=1}^{n-1}\Delta_k\Big)\cup\Delta_{n+1},
$$
where
$$
\Delta_k:=\Big\lbrace\lambda\in\sigma_{re}(D_k):\  \sum\limits_{s=k+1}^{n}\alpha(D_s-\lambda)<\infty\Big\rbrace,\ 1\leq k\leq n-1.
$$
$$
\Delta_{n+1}=\lbrace\lambda\in\mathds{C}:\ \sum\limits_{s=1}^n\alpha(D_s-\lambda)<\sum\limits_{s=1}^n\beta(D_s-\lambda)\rbrace.
$$
\end{Posledica}

\subsection{When does $\sigma_*(T_n^d(A))=\bigcup\limits_{s=1}^k\sigma_*(D_s)$ hold for arbitrary $A\in\mathcal{B}_n$, $\sigma_*\in\lbrace \sigma_{lw},\sigma_{rw}\rbrace$?}
In this subsection we provide an answer to Question 3 stated in Section 2. This question has been solved for case of infinite dimensional separable Hilbert spaces by Wu and Huang in \cite{WU2}. Again, it turns out that their results stated in Section 3 of \cite{WU2} hold, word by word, in arbitrary Banach spaces, with the only distinction that in every claim we have to assume validity of the complements condition for all $D_i$'s. We invite readers to copy the statements of Theorems 3.1, 3.6 from \cite{WU2} and their Corollaries, by including the complements condition assumption and substituting Hilbert spaces $\mathcal{H}_1$,..., $\mathcal{H}_n$ by Banach spaces $X_1,...,X_n$. As of special interest, we recommend studying Remark 3.5 of \cite{WU2} in a connection with an inte\-re\-sting notion of single valued extension property of Finch (\cite{FINC}). However, for the convenience of a reader, we decided to provide a couple of mentioned statements for illustration, and so we give analogues of Corollaries 3.3 and 3.8. Their proofs can be found in \cite{WU2}, Section 3.
\begin{Posledica}
Let $D_1\in\mathcal{B}(X_1),\ D_2\in\mathcal{B}(X_2),...,D_n\in\mathcal{B}(X_n)$ with all $D_i,\ i=1,...,n $ satisfying the complements condition. Then
$$
\sigma_{lw}(T_n^d(A))=\bigcup\limits_{k=1}^n\sigma_{lw}(D_k)
$$
holds for every $A\in\mathcal{B}_n$ if and only if
$$
\Delta_1\cup\Delta_2=\emptyset
$$
where $\Delta_k$, $k=1,2$ are defined as follows
$$
\Delta_1:=\bigcup\limits_{k=2}^n\lbrace\lambda\in\bigcap\limits_{s=1}^{k-1}\rho_{le}(D_s)\cap[\bigcup\limits_{s=1,\ s\neq k-1}^n\sigma_{lw}(D_s)]:\ \beta(D_{k-1}-\lambda)=\infty\rbrace,
$$
$$
\begin{aligned}
\Delta_2:=\bigcup\limits_{k=1}^n\lbrace\lambda\in\bigcap\limits_{s=1}^{n}\rho_{le}(D_s):\ \sum\limits_{s=1}^n\alpha(D_s-\lambda)\leq\sum\limits_{s=1}^n\beta(D_s-\lambda),\\ \alpha(D_k-\lambda)>\beta(D_k-\lambda)\rbrace.
\end{aligned}
$$
\end{Posledica}
\begin{Posledica}
Let $D_1\in\mathcal{B}(X_1),\ D_2\in\mathcal{B}(X_2),...,D_n\in\mathcal{B}(X_n)$ with all $D_i,\ i=1,...,n $ satisfying the complements condition. Then
$$
\sigma_{rw}(T_n^d(A))=\bigcup\limits_{k=1}^n\sigma_{rw}(D_k)
$$
holds for every $A\in\mathcal{B}_n$ if and only if
$$
\Delta_1\bigcup\limits\Delta_2=\emptyset
$$
where $\Delta_k$, $k=1,2$ are defined as follows
$$
\Delta_1:=\bigcup\limits_{k=1}^{n-1}\lbrace\lambda\in\bigcap\limits_{s=k+1}^{n}\rho_{re}(D_s)\cap[\bigcup\limits_{s=1,\ s\neq k+1}^n\sigma_{rw}(D_s)]:\ \alpha(D_{k+1}-\lambda)=\infty\rbrace,
$$
$$
\begin{aligned}
\Delta_2:=\bigcup\limits_{k=1}^n\lbrace\lambda\in\bigcap\limits_{s=1}^{n}\rho_{re}(D_s):\ \sum\limits_{s=1}^n\beta(D_s-\lambda)\leq\sum\limits_{s=1}^n\alpha(D_s-\lambda),\\ \beta(D_k-\lambda)>\alpha(D_k-\lambda)\rbrace.
\end{aligned}
$$
\end{Posledica}

\noindent\textbf{Acknowledgments}\\[3mm]
\hspace*{6mm}This work was supported by the Ministry of Education, Science and Technological Development of the Republic of Serbia under Grant No. 451-03-9/2021-14/ 200125\\[1mm]
\hspace*{6mm}I wish to express my gratitude to professor Dragan S. Djordjevi\'{c} for introduction to this topic and for useful comments that greatly improved the form of this paper.  \\[3mm]


\begin{thebibliography}{W}

\bibitem{APOSTOL} Apostol C., {\it The reduced minimum modules}, Michigan Math. J. \textbf{32}, 279 - 294 (1985)

\bibitem{NEOGRANICENI3} Bai Q. \& Chen A. \& Huang J., {\it Some spectra properties of unbounded 2×2 upper triangular operator matrices,} Ann. Funct. Anal. \textbf{10} (2019), no. 3, 412–424

\bibitem{CAO} Cao X. H. \& Guo M. Zh. \& Meng B., {\it Semi-Fredholm spectrum and Weyl's theorem for operator matrices,} Acta Math. Sin. (Engl. Ser.) \textbf{22} (2006), no. 1, 169–178

\bibitem{CAO2}  Cao X. \& Meng B., {\it Essential approximate point spectra and Weyl's theorem for operator matrices}, J. Math. Anal. Appl. \textbf{304} (2005), no. 2, 759–771

\bibitem{CARADUS}  Caradus S. R. \& Pfaffenberger W. E. \& Yood, B., {\it Calkin algebras and algebras of operators on Banach spaces,} Lecture Notes in Pure and Applied Mathematics, Vol. 9. Marcel Dekker, Inc., New York, 1974. viii+146 pp.

\bibitem{SUNTHEORY}  Djordjevi\'c D. S., {\it Perturbations of spectra of operator matrices}, J. Oper. Theory. \textbf{48}(3), 467-486 (2002)

\bibitem{DJR}   Djordjevi\'c D. S. \& Rako\v cevi\'c V., {\it  Lectures on generalized inverses}, University of Ni\v s,
Faculty of Sciences and Mathematics, Ni\v s, 2008

\bibitem{DU} Du H. K. \& Jin P., {\it Perturbation of spectrums of 2×2 operator matrices}, Proc. Amer. Math. Soc. \textbf{121} (1994), no. 3, 761–766

\bibitem{DUNFORD}  Dunford, N. \& Schwartz, J. T., {\it Linear operators. Part I. General theory. With the assistance of William G. Bade and Robert G. Bartle.} Reprint of the 1958 original. Wiley Classics Library. A Wiley-Interscience Publication. John Wiley \& Sons, Inc., New York, 1988

\bibitem{FINC} Finch J. K., {\it The single valued extension property on a Banach space}. Pacific J. Math. \textbf{58} (1975), no. 1, 61–69

\bibitem{HAN} Han J. K. \& Lee H. Y. \& Lee W. Y., {\it Invertible completions of 2×2 upper triangular operator matrices}, Proc. Amer. Math. Soc. \textbf{128} (2000), no. 1, 119–123

\bibitem{HUANG} Huang J. \& Wu X. \& Chen A., {\it The point spectrum, residual spectrum and continuous spectrum of upper-triangular operator matrices with given diagonal entries}, Mediterr. J. Math. \textbf{13} (2016), no. 5, 3091–3100

\bibitem{KATO} Kato T., {\it Perturbation theory for linear operators.} Reprint of the 1980 edition. Classics in Mathematics. Springer-Verlag, Berlin, 1995

\bibitem{LI} Lee W. Y., {\it Weyl spectra of operator matrices}, Proc. Amer. Math. Soc. \textbf{129} (2001), no. 1, 131–138

\bibitem{KINEZI} Li Y. \& Sun X. H. \& Du H. K., {\it The intersection of left (right) spectra of 2×2 upper triangular operator matrices}. Linear Algebra Appl. \textbf{418} (2006), no. 1, 112–121

\bibitem{SUN}  Li Y. \& Sun X.H. \& Du H.K. {\it Intersections of the left and right essential spectra of 2×2 upper triangular operator matrices}. Bull. London Math. Soc. \textbf{36} (2004), no. 6, 811–819

\bibitem{LIandDU} Li Y. \& Du H., {\it The intersection of essential approximate point spectra of operator matrices}. J. Math. Anal. Appl. \textbf{323} (2006), no. 2, 1171–1183

\bibitem{NEOGRANICENI2} Liu J. \& Huang J. \& Chen A., {\it Spectral inclusion by the quadratic numerical range of 2×2 operator matrices with unbounded entries.} Filomat \textbf{34} (2020), no. 4, 1283–1293

\bibitem{MULLER} M\"{u}ller V., {\it Spectral Theory of Linear Operators and Spectral Systems in Banach Algebras}, Birkh\"{a}user - Verlag, Basel (2003)

\bibitem{NEOGRANICENI1} Rasulov T. \& Tretter C. {\it Spectral inclusion for unbounded diagonally dominant n×n operator matrices}. Rocky Mountain J. Math. \textbf{48} (2018), no. 1, 279–324

\bibitem{SETER} Schechter M., {\it Principles of Functional Analysis}, Academic Press, New York, 1971.

\bibitem{WU} Wu X. \& Huang J., {\it Essential spectrum of upper triangular operator matrices}, Ann. Funct. Anal, \textbf{11} (2020), no. 3, 780–798

\bibitem{WU2} Wu X. \& Huang J., {\it Weyl spectrum of upper triangular operator matrices}, Acta Math. Sin. \textbf{36}(7), 783-796 (2020)

\bibitem{ZGUITTI} Zguitti H., {\it A note on Drazin invertibility for upper triangular block operators}, Mediterr. J. Math. \textbf{10} (2013), no. 3, 1497–1507

\bibitem{ZHANG} Zhang S. \& Wu Z., {\it Characterizations of perturbations of spectra of 2×2 upper triangular operator matrices.} J. Math. Anal. Appl. \textbf{392} (2012), no. 2, 103–110

\bibitem{ZANA} \v{Z}ivkovi\'{c} Zlatanovi\'{c} S. \v{C} \& Rako\v{c}evi\'{c} V. \& Djordjevi\'{c} D. S., {\it Fredholm theory}, to appear
\end{thebibliography}
\end{document}